\documentclass[12pt]{amsart}
\usepackage{amsbsy}
\begin{document}

\newtheorem{theorem}{Theorem}
\newtheorem{proposition}{Proposition}
\newtheorem{lemma}{Lemma}
\newtheorem{corollary}{Corollary}
\newtheorem{definition}{Definition}
\newtheorem{remark}{Remark}
\numberwithin{equation}{section}
\numberwithin{theorem}{section}
\numberwithin{proposition}{section}
\numberwithin{lemma}{section}
\numberwithin{corollary}{section}
\numberwithin{definition}{section}
\numberwithin{remark}{section}
\newcommand{\ren}{{\mathbb R}^N}
\newcommand{\re}{{\mathbb R}}
\newcommand{\n}{\nabla}
\newcommand{\iy}{\infty}
\newcommand{\pa}{\partial}
\newcommand{\fp}{\noindent}
\newcommand{\ms}{\medskip\vskip-.1cm}
\newcommand{\mpb}{\medskip}

\newcommand{\tex}{\textstyle}
\renewcommand{\a}{\alpha}
\renewcommand{\b}{\beta}
\newcommand{\g}{\gamma}
\newcommand{\G}{\Gamma}
\renewcommand{\d}{\delta}
\newcommand{\D}{\Delta}
\newcommand{\e}{\varepsilon}
\renewcommand{\l}{\lambda}
\renewcommand{\o}{\omega}
\renewcommand{\O}{\Omega}
\newcommand{\s}{\sigma}
\renewcommand{\t}{\tau}
\newcommand{\z}{z}
\newcommand{\wx}{\widetilde x}
\newcommand{\wt}{\widetilde t}
\newcommand{\noi}{\noindent}
\newcommand{\AAA}{{\mathbf   A}}
\newcommand{\ee}{\varepsilon}
\newcommand{\var}{\varphi}
\newcommand{\BB}{{\mathbf  B}}
\newcommand{\CC}{{\mathbf  C}}
\newcommand{\II}{{\mathrm{Im}}\,}
\newcommand{\RR}{{\mathrm{Re}}\,}
\newcommand{\eee}{{\mathrm  e}}
\newcommand{\uu}{{\bf u}}
\newcommand{\xx}{{\bf x}}
\newcommand{\yy}{{\bf y}}
\newcommand{\zz}{{\bf z}}
\newcommand{\aaa}{{\bf a}}
\newcommand{\cc}{{\bf c}}
\newcommand{\jj}{{\bf j}}
\newcommand{\ggg}{{\bf g}}
\newcommand{\UU}{{\bf U}}
\newcommand{\YY}{{\bf Y}}
\newcommand{\HH}{{\bf H}}
\newcommand{\GGG}{{\bf G}}
\newcommand{\VV}{{\bf V}}
\newcommand{\ww}{{\bf w}}
\newcommand{\vv}{{\bf v}}
\newcommand{\hh}{{\bf h}}
\newcommand{\di}{{\rm div}\,}
\newcommand{\LL}{L^2_\rho(\ren)}
\newcommand{\LLL}{L^2_{\rho^*}(\ren)}
\newcommand{\inA}{\quad \mbox{in} \quad \ren \times \re_+}
\newcommand{\inB}{\quad \mbox{in} \quad}
\newcommand{\atA}{\quad \mbox{at} \quad}
\newcommand{\onA}{\quad \mbox{on} \quad}
\newcommand{\inC}{\quad \mbox{in} \quad \re \times \re_+}
\newcommand{\inD}{\quad \mbox{in} \quad \re}
\newcommand{\forA}{\quad \mbox{for} \quad}
\newcommand{\whereA}{,\quad \mbox{where} \quad}
\newcommand{\asA}{\quad \mbox{as} \quad}
\newcommand{\andA}{\quad \mbox{and} \quad}
\newcommand{\withA}{,\quad \mbox{with} \quad}
\newcommand{\orA}{,\quad \mbox{or} \quad}
\newcommand{\ef}{\eqref}
\newcommand{\ssk}{\smallskip}
\newcommand{\LongA}{\quad \Longrightarrow \quad}
\def\com#1{\fbox{\parbox{6in}{\texttt{#1}}}}
\def\N{{\mathbb N}}
\def\A{{\cal A}}
\newcommand{\de}{\,d}
\newcommand{\eps}{\varepsilon}
\newcommand{\beq}{\begin{equation}}
\newcommand{\eeq}{\end{equation}}
\newcommand{\spt}{{\mbox spt}}
\newcommand{\ind}{{\mbox ind}}
\newcommand{\supp}{{\mbox supp}}
\newcommand{\dip}{\displaystyle}
\newcommand{\prt}{\partial}
\newcommand{\ii}{{\rm i}\,}
\newcommand{\Dm}{(-\D)^m}

\title
{\bf Refined scattering and  Hermitian spectral theory for linear
higher-order
 Schr\"odinger equations}

\author {V.A.~Galaktionov and I.V.~Kamotski}

\address{Department of Mathematical Sciences, University of Bath,
 Bath BA2 7AY, UK}
\email{vag@maths.bath.ac.uk}

\address{Department of Mathematical Sciences, University of Bath,
 Bath BA2 7AY, UK}
\email{ivk20@maths.bath.ac.uk}



\date{\today}

\begin{abstract}

The Cauchy problem for a linear $2m$th-order Schr\"odinger
equation
 \beq
 \label{0.1}
 u_t= -{\rm i} \,(-\D)^m u \inB \ren \times \re_+, \quad
 u|_{t=0}=u_0; \quad m \ge 1 \,\,\, \mbox{is an integer},
 \eeq
 is studied, for initial data $u_0$ in the weighted space $L^2_{\rho^*}(\ren)$, with
 $\rho^*(x)={\mathrm e}^{|x|^\a}$ and $\a= \frac {2m}{2m-1}\in (1,2]$.
 The following {\bf five} problems are studied:

 \noi {\bf (I)} A sharp asymptotic behaviour of solutions as $t \to +\iy$  is
 governed by a discrete spectrum and a countable set $\Phi$ of
the eigenfunctions of the linear rescaled operator
 $$
  \tex{
  \BB = -\ii \Dm + \frac 1{2m} \, y \cdot \n + \frac {N}{2m}\, I,
  \,\,\,\mbox{with the spectrum}\,\,\,
   \s(\BB)= \big\{ \l_\b= - \frac{|\b|}{2m}, \, |\b| \ge
   0\big\}.
  }
  $$

\noi {\bf (II)} Finite-time blow-up local structures of  nodal
sets of solutions as $t \to 0^-$ and formation of ``multiple
zeros" are
 described by the eigenfunctions being   {\em generalized Hermite polynomials}, of the
 ``adjoint" operator
 $$
  \tex{
  \BB^* = -\ii \Dm - \frac 1{2m} \, y \cdot \n,
   \,\,\,\mbox{with the same spectrum}\,\,\,
  \s(\BB^*)=\s(\BB).
  }
  $$

Applications of these spectral results also  include: {\bf (III)}
a unique continuation theorem,  and {\bf (IV)} boundary
characteristic point regularity issues.

Some applications
are discussed for more general linear PDEs and for the {nonlinear
Schr\"odinger equations}
in the focusing (``$+$") and defocusing (``$-$") cases
 $$
 u_t= -{\rm i} \,(-\D)^m u \pm \ii |u|^{p-1}u \inB \ren \times \re_+,
 \quad \mbox{where} \quad
  p>1,
  $$
  as well as for {\bf (V)} the {quasilinear} Schr\"odinger equation of
   a ``porous medium type"
   $$
   u_t=-\ii (-\D)^m(|u|^n u) \inA \whereA n>0.
   $$
   For the latter one,
 the main idea towards countable families of {\em
 nonlinear eigenfunctions} is to perform a homotopic path $n \to
 0^+$ and to use  spectral theory of the pair $\{\BB,\BB^*\}$.

\end{abstract}

\maketitle

\section{Introduction: duality of  global and blow-up scalings, Hermitian spectral
 theory, and refined scattering}
\label{S1}

\subsection{Basic Shr\"odinger equations and key references}

Consider the linear $2m$th-order Schr\"odinger equation (the
LSE--$2m$), with  any integer $m \ge 1$,
 \beq
 \label{1}
 u_t= -{\rm i}\, (-\D)^m u \inB \ren \times \re_+, \quad u|_{t=0}=u_0,
  \eeq
   where $\D$ is the Laplace operator in $\ren$,
  for initial data $u_0$ in some weighted $L^2$-space, to be
  introduced.
   Here $m=1$
  corresponds to the classic Schr\"odinger equation
 \beq
 \label{2}
 \ii u_t= -\D u \inB \ren \times \re_+,
  \eeq
   which very actively entered general PDE theory from Quantum Mechanics in the 1920s.

It is not  an exaggeration to say that, nowadays, linear and
nonlinear Schr\"odinger type equations are the most popular PDE
models of modern mathematics among other types of equations. In
Appendix A, we  present a ``mathematical evidence" for that by
using simple data from the {\tt MathSciNet}.
 It not possible to express how deep is mathematical theory
developed for models such as \ef{2}, \ef{1}, and related
semilinear ones. We refer to well-known monographs \cite{SulMon,
Caz03}, which cover classes of both linear and nonlinear PDEs.

 Concerning the results that are more closely related to the subject of this paper, we note
 that scattering  $L^2$ and $L^{q,r}_{x,t}$ theory for \ef{2} has been fully developed
in the works by Stein, Tomas, Segal, Strichartz   in the 1970s
 with later further involved  estimates in more general spaces by Ginibre and
 Velo, Yajima, Cazenave and Weissler, Montgomery-Smith, Keel, Tao,
 and many others; see \cite{Keel98} and \cite{Vil07} for references concerning these,
  as well as  optimal  $L^{q,r}_{x,t}$  estimates
 for the non-homogeneous equation
\beq
 \label{2F}
 \ii u_t= -\D u + F(x,t) \inB \ren \times \re,
  \eeq
  as well as more recent papers \cite{Koh11, Rogers09, Zhai09}.

 The $2m$th-order counterpart \ef{1} was also under scrutiny for
 a long period. We refer to  Ablowitz--Segur's monograph
    \cite{Abl81}, Ivano--Kosevich \cite{Ivano83}, Turitsyn \cite{Tur85},
     Karpman \cite{Karp96}, and Karpman--Shagalov \cite{Karp00} for physical,
     symmetry, and other backgrounds
     of higher-order
    Schr\"odinger-type semilinear models (see also \cite{Zhu07} for extra motivations
     from nonlinear optics), \cite{Pech79} for first existence and uniqueness results,
     and more recent papers
    \cite{Bar10, Bar10Ring, Cui07, Guo10, Hao05, Hao07, Miao09, Paus09, Paus10, Zhu10} as an account for further applied and
rigorous research, as well as other earlier key references and
survey, in this fundamental area of modern PDE theory.

\subsection{Discrete real spectra, ``Hermitian spectral history", and our motivation}

Actually, the developed here refined scattering theory is rather
general, so our methods can be applied not only to the
Schr\"odinger equations such as \ef{1}, but also to practically
any linear evolution PDEs with constant or perturbed smooth
coefficients and classic solutions. Formulating the approach
rather loosely, we claim that, for \ef{1}, the most principal part
is played by the spectral theory for the following rescaled
operator:
 \beq
 \label{RO1}
  \tex{
 \BB^*= - \ii (-\D_y)^m - \frac 1{2m}\, y \cdot \n_y
 \whereA y= \frac x{(-t)^{1/2m}} \forA t<0
 }
  \eeq
  is the corresponding  {\em Sturmian}  {\em blow-up backward
  spatial variable}  at the focusing point $(x,t)=(0,0^-)$.
In 1836, C. Sturm  used the backward
  variable $y=\frac x{\sqrt{-t}}$, $t \to 0^-$, for the heat equation with a potential,
  \beq
  \label{HE11}
 u_t=u_{xx} + q(x) u,
 \eeq
 in his seminal paper
  \cite{St} [The pioneering work \cite{St} was practically fully forgotten for about
  150 years!--practically until the 1980s,--plausibly the most amazing and striking such an
  example in the whole history of mathematics ever], where he formulated his
  two fundamental  theorems on zeros sets of
  solutions $u=u(x,t)$ of \ef{HE11}. This remarkable history, with many extensions about, can be found in
   \cite[Ch.~1]{GalGeom} with
  precise statements of Sturm's results of the 1836 written in his original notations.

   The operator \ef{RO1} and its adjoint below are
   then respectively
    defined in  weighted $L^2$-spaces, with a
special ``radiation conditions at infinity" to be specified.
 More than half of the paper is devoted to the study of \ef{RO1} and its
``adjoint" operator
 \beq
 \label{BBN}
  \tex{
 \BB=- \ii (-\D_y)^m + \frac 1{2m}\, y \cdot \n_y+ \frac N{2m} \,
 I , \quad
 y= \frac x{t^{1/2m}}, \,\,\, t>0\,\,\, (\mbox{the {\em forward variable}}).
 }
 \eeq

Indeed,  \ef{RO1} is a perturbation of the original one in \ef{1}.
Though the perturbation is of the first order, the coefficient $y$
therein is unbounded as $y \to + \iy$, so this  changes the
natural space $L^2(\ren)$ and moves the operator into an
essentially weighted metric. An amazing property of \ef{RO1} is
that, being properly defined, it has the discrete spectrum
 \beq
  \label{RO2}
  \tex{
   \s(\BB^*)=\big\{ \l_\b=- \frac{|\b|}{2m},
   \,\,|\b|=0,1,2,...\big\}
   \quad (\mbox{$\b$ is a multiindex in $\ren$}),
    }
    \eeq
and all the eigenfunctions $\{\psi_\b^*(y)\}$ are finite
polynomials ({\em generalized Hermite} ones).


These properties directly match the classic results for the {\em
heat equation} for $m=1$:
 \beq
 \label{RO3}
  \begin{matrix}
 u_t=u_{xx} \LongA \BB^*= D^2_y - \frac 12 \, y D_y,
  \,\,\, \mbox{with} \,\,\, \s(\BB^*)=\big\{ \l_\b=- \frac{l}{2},
  \,\,l \ge 0\big\}, \qquad\quad
  \ssk\ssk\\
   \mbox{so}
  \,\,\, (\BB^*)^*=\BB^* \inB L^2_{\hat \rho^*}(\re), \,\,\,\hat \rho^*(y)={\mathrm e}^{-\frac
  {y^2}4},\,\,\,
\BB=D^2_y + \frac 12 \, y D_y+ \frac 12\, I,
  \qquad\qquad
 \end{matrix}
  \eeq
  where $\BB$ is defined in the adjoint space $L^2_{\hat \rho}$, with
  $\hat \rho=\frac 1{\hat \rho^*}$, and $(\BB)^*=\BB^*$ in the dual metric of $L^2$, etc.
In  modern language, for \ef{RO3}, the spectrum \ef{RO2} with
$m=1$ and the classic Hermite polynomials (introduced in detail
about 1870) as eigenfunctions were already constructed by C.~Sturm
in 1836 \cite{St}.
 As we  mentioned already, this  led Strum to formulate his two fundamental theorems on the
structure of multiple zeros of solutions of parabolic equations
and on nonincrease in time of the zero number (or sign changes of
solutions);
 we refer again to  \cite[Ch.~1]{GalGeom} for a full
history and key further references and extensions. The operator
$\BB^*$ in \ef{RO3}, admitting a natural $N$-dimensional extension
similar to \ef{RO1}, remains one of the   key objectives  in
general theory of linear self-adjoint operators; see
Birman--Solomjak's monograph \cite{BS}.


The spectral results for \ef{RO3} and their consequences for the
asymptotic behaviour for second-order parabolic equations are
classic and well-known since the 1830s, with further extensions as
orthonormal polynomial families by Hermite himself from the 1870s.


 However, and this looks like a truly amazing fact, a direct extension
of such classic results
 to other classes of
PDEs took a lot of time. For instance, similar spectral theory for
the 1D {\em bi-harmonic equation} (cf. with the one in the  first
line in \ef{RO3})
 \beq
 \label{RO4}
 \tex{
 u_t=-u_{xxxx}\,\, \Longrightarrow \,\, \BB^*= -D^4_y - \frac 14 \, y D_y,
  \,\,\mbox{with} \,\, \s(\BB^*)=\big\{ \l_\b=- \frac{l}{4},
   \,\,l=0,1,2,...\big\},
 }
  \eeq
  etc.,
was developed in 2004 \cite{Eg4}, i.e.,  168 years later after
Sturm's pioneering discovery for $m=1$ in 1836! As a certain (but
seems not that convincing) excuse, note that the operator $\BB^*$
in \ef{RO4} {\em is not self-adjoint} in no weighted space, though
keeps having discrete real spectrum, polynomial eigenfunctions
(naturally called {\em generalized Hermite ones}), and a number of
other nice and typical from self-adjoint theory properties.


Therefore, our goal is to  show that similar issues  remain true
for our rescaled Schr\"odinger operators \ef{RO1} and  \ef{BBN}
(so we are talking about a {\em spectral pair} $\{\BB,\BB^*\}$ of
non-self-adjoint ones), which clearly have  analogous structures,
though the mathematics becomes essentially more involved than for
\ef{RO4}, to say nothing of the well-studied self-adjoint case
\ef{RO3}.




\subsection{Layout of the paper: duality of global  and blow-up  asymptotics}

  In Section \ref{S2}, we describe some
properties of the {\em fundamental solution} of \ef{1} given by
 \beq
 \label{Fund}
  \tex{
  b(x,t)= t^{-\frac N{2m}} \, F(y), \quad y = \frac x{t^{1/{2m}}}
 \LongA \BB F=0 \inB \ren
  .
  }
   \eeq
   Next sections are mainly  devoted to the following two asymptotic problems
    for \ef{1}:

 \ssk

\noi{\bf Application I:} {\sc Global asymptotics as $t \to +\iy$,
Sections \ref{S3}
 and \ref{S5Compl}}. The
asymptotic behaviour as $t \to +\iy$ of solutions is  governed by
the eigenfunctions of \ef{BBN}:
 \beq
 \label{3}
   \tex{
  \BB = -\ii \Dm + \frac 1{2m} \, y \cdot \n + \frac {N}{2m}\, I,
 \,\,\, \mbox{with} \,\,\,
    \s(\BB)= \big\{ \l_\b= - \frac{|\b|}{2m}, \, |\b| \ge
   0\big\},
  }
  \eeq
  which demands a proper definition of its domain by a spectral decomposition (Section \ref{S3})
  and by a traditional spectral theory involving careful using poles of the resolvent,
  non-classic ``radiation conditions" at infinity  (see an
  alternative approach in \cite{GK2mClass}), etc.
 More precisely, we establish that the discrete spectrum and  the eigenfunction set for the
 operator \ef{3} describe  all the possible asymptotics as $t \to
 +\iy$ of solutions of \ef{1} for any data $u_0 \in
 L^2_{\rho^*}(\ren)$.
  The
exponential  weights (recall the change $\rho \mapsto \rho^*$
relative to \ef{RO3} done by some clear and natural reasons)
 \beq
 \label{weights}
  \tex{
 \rho^*(y)={\mathrm e}^{|y|^\a} \withA \a=
\frac{2m}{2m-1}, \andA \rho(y)= \frac 1{\rho^*(y)}={\mathrm
e}^{-|y|^\a},
 }
 \eeq
  are properly
introduced in Section \ref{S3}.
 It is curious that, even in the classic case $m=1$, i.e., for
 \ef{2}, we were not able to find any essential traces of such a full refined
 scattering theory (except some  particular results often admitting not-that-clear interpretation)
 and corresponding  spectral
 properties in the vast existing scattering literature\footnote{The authors
   do not still believe that
 optimal and sharp large-time ($t \to \iy$) asymptotic  theory
 for the classic LSE \ef{2} in $\ren$ has not been  developed  in
  {\em full details} (some particular results have been indeed known seems) since this
 PDE burst   into  quantum mechanics and mathematical
 physics in the 1920s.
 As was mentioned, a similar (and even stronger) asymptotic theory for the 1D heat equation
 $u_t=\D u$ (cf. \ef{1} for $m=1$) is known since Sturm's analysis
 \cite{St} of formation of multiple zeros of solutions
 (we will do this for \ef{1} in Section \ref{S6}) obtained
 in 1836!
  The authors will be very pleased to get rid of
 such a quite
 surprising
 delusion, but also would be naturally satisfied to know that this, though again unbelievably,
  is done in full details  for the first time
 in the present paper.}.


 \ssk

  The classic real analogy and forerunner  of \ef{3} is
 the self-adjoint operator for $m=1$ (we apologize for the
 necessary change of the weights here, $\rho \mapsto \rho^*$; cf.
 \ef{weights})
 \beq
 \label{HH1}
 \tex{
\BB = \D + \frac 1{2} \, y \cdot \n + \frac {N}{2}\, I \equiv
\frac 1 \rho \, \n \cdot (\rho \n)+ \frac N2 \, I
 \inB \LL, \quad \mbox{with} \quad \rho={\mathrm e}^{
 \frac{|y|^2}4}.
  }
  \eeq
 As we have pointed out already, its real discrete spectrum $\s(\BB)=\{ - \frac l2, \, l \ge 0\}$
   and eigenfunctions as {\em Hermite
  polynomials}, multiplied by the Gaussian, are known from, at least, the 1830s, and are
  associated with the names of Sturm  and Hermite;
  see \cite[\S~1.2]{GalGeom} for more history and original Sturm's calculations,
   and \cite[p.~48]{BS} for a fuller
  account of applications of
   these separable polynomials in self-adjoint linear operator theory.

Thus, we are obliged  here to develop  Hermitian-like spectral
theory for the $2m$th-order rescaled Schr\"odinger operator
\ef{3}, and this is an unavoidable task if we want to reach an
optimal classification of large time behaviour for the
non-stationary LSE \ef{1}.


 \ssk

\noi {\bf Application II:} {\sc Blow-up asymptotics as $t \to
0^-$, Sections \ref{S4},
 \ref{S6}, Appendices  B and D}.
Alternatively, for data $u_0 \in \LLL$, using
 blow-up
scaling
 at a finite point as $x \to 0$ and $t \to  T^-=0^-$, we
 show that this
behaviour of solutions and local structure of their nodal sets are
 described by the eigenfunctions of the linear operator \ef{RO1}, which is
  ``adjoint" to $\BB$ in a  sense (but not in the standard
 dual $L^2$-metric; an {\em indefinite metric} should be involved, which will be carefully explained).
 The discrete real spectrum \ef{RO2} remains the same as in \ef{3}.

 A key point is that  the eigenfunctions of $\BB^*$ are {\em
 generalized
 Hermite   polynomials} $\{\psi_\b^*(y)\}$, so that the nodal sets of solutions
 of \ef{1} are locally governed by zero surfaces generated by
 these polynomials {\em only} (there exists a countable, complete,
 and closed set of those).
 {\bf Application III:} This allows us to state in Section
 \ref{S6} a sharp {\em uniqueness continuation} theorem for
 \ef{1}: if  $\II u(x,t)$ (or $\RR u(x,t)$) has, roughly speaking,
  \beq
  \label{5}
   \begin{matrix}
  \mbox{a rescaled nodal set component not decomposable into a finite
  combination}\qquad
  \ssk\\
  \mbox{from polynomial surfaces $\{\II \psi_\b^*(y)=0, \,
  |\b| \ge 0\}$, then $u(x,t) \equiv 0$.}\qquad
   \end{matrix}
   \eeq
Some further applications of these spectral results are also
 discussed.

\ssk

   In Sections \ref{S5Compl} and \ref{S6}, possible applications of
    Hermitian spectral theory are studied for more general linear PDEs and for the
 $2m$th-order  nonlinear Schr\"odinger equation (the NLSE--$2m$)
   \beq
   \label{6}
   u_t= -\ii (-\D)^m u \pm   \ii |u|^{p-1} u \inB \ren \times \re_+ \whereA
   p>1
    \eeq
 and the sign ``$+$" corresponds to the focusing (blow-up) model,
 while ``$-$" gives a defocusing one. See
 \cite{Ken06, Mer04, Mer05, Plan07, Rap06, Vis107} as a guide concerning the modern research
 of both semilinear PDEs \ef{6}.


 {\bf Application IV:}
In Appendix B,
we show how to apply the spectral results to the classic problem
of the {\em regularity of a boundary characteristic point} for the
linear Schr\"odinger equation. It seems that, in the present
setting of rather arbitrary ``backward paraboloids" at
characteristic points, such issues were not addressed in the
existing  literature.

 \ssk

\noi(ii) {\bf Application V:} {\sc Quasilinear Schr\"odinger
equation, Appendix C}. This is a most ``risky" application of our
refined scattering spectral theory for the spectral pair
$\{\BB,\BB^*\}$, so we also put this into an
appendix\footnote{This and some other non-entirely-rigorous parts
of our applications well resonate with two
Tao's  comments in his ``What is good mathematics?":  \noi (i)
 ``...
mathematical rigour, while highly important, is only  one
component of what determines a quality piece of
 mathematics" \cite[p.~624]{TTaoMath07}, and
\noi (ii) ``... we should also be aware of any possible
  larger context that one's results could be placed in, as this may lead to
  the greatest long-term benefit for the result, for the field, and for mathematics
  as a whole." \cite[p.~633, Tao's last sentence therein]{TTaoMath07}.
  [In author's opinion, Tao's (ii) statement could definitely
  serve as a most impressive characterization of a ``good
  mathematics" among dozens of others presented in his paper.]
 The latter one (ii) somehow helps us to justify  including some ``risky
 applications", which we are not aware how to prove, and even are
 not sure whether such results can be proved in any remote  future;
 though the authors believe that this stuff should be revealed
 for the Readers, which might be interested in nowadays or will be later on.
 }. Thus, exhibiting a certain necessary bravery, we develop some
basics of a ``nonlinear eigenfunction theory" for a {\em
quasilinear Schr\"odinger equation} (the QLSE)
 of the form
 \beq
 \label{QQ1}
  u_t= - \ii (-\D)^m (|u|^n u) \inB \ren \times \re_+ \whereA n>0
   \eeq
   is a fixed parameter. The applications and some history/references concerning such rather unusual
    quasilinear PDEs are
   explained therein. Here, we intend to reconstruct a proper
   connection between linear and ``nonlinear" spectral theory by
   performing a continuity homotopic path $n \to 0^+$, which
   establishes a  link between  ``nonlinear eigenvalue problems"
for \ef{QQ1} and the linear one developed for \ef{1}. As a result,
we predict existence of a countable family of the so-called
$n$-branches of solutions, which are originated at $n=0$ from
eigenfunctions on
 the corresponding eigenspaces for the linear spectral pair
 $\{\BB,\BB^*\}$.

\ssk

We strongly believe that the results of Hermitian spectral theory
developed can be useful for attacking a number of open problems
concerning blow-up behaviour for \ef{6}, \ef{QQ1}, and others. We
plan to explain this
 in a forthcoming paper. Meantime, we
just comment on that standard blow-up rescaling ({\em q.v.}
\ef{RR1}) leads to the adjoint operator \ef{RO1} as the
linearization, so the generalized Hermite polynomial
eigenfunctions of $\BB^*$ can be key for understanding this
intriguing rescaled blow-up dynamics.

 \section{Fundamental solution and
 the convolution}
 \label{S2}

\subsection{Fundamental solution and its first properties}

In constructing fundamental solutions and corresponding
convolutions for \ef{1}, one can use the fact that formally
changing the independent time variable $\ii t \mapsto t$ yields
the standard poly-harmonic PDE:
 \beq
 \label{2.1}
 \ii t \mapsto t \LongA u_t= - \Dm u \quad \mbox{in} \quad \ren
 \times \re_+.
  \eeq
  This creates  an artificial complex (imaginary) time axis, and nevertheless will help us to
  restore various spectral properties and other functional details
 related to \ef{1}. Of course, unlike the real parabolic case \cite{Eg4},
  the change in \ef{2.1} implies a well known highly oscillatory
  properties of the fundamental and other
 solutions of the Schr\"odinger equation \ef{1}, that are not
 available for its real parabolic counterpart in \ef{2.1}.

Thus, by classic PDE theory, given proper initial data $u_0(x)$,
the unique solution of the Cauchy problem for (\ref{2.1}) is given
by
 \beq
 \label{CP11}
  u(t)= b(t) * u_0 \forA t>0,
   \eeq
   where $b(x,t)$ is the fundamental solution (\ref{Fund}) of the operator in
   (\ref{2.1}).
Substituting $b(x,t)$ into (\ref{1}), one obtains the {\em
rescaled kernel}
 $F(y)$ as a unique solution of a complex linear ordinary differential equation (ODE) which is the
radial restriction of a linear PDE system,
\begin{equation}
\label{ODEf}
 \tex{
 {\mathbf B}F \equiv - \ii (-\Delta _y)^m F + \frac
1{2m}\, y \cdot \nabla_y  F + \frac N{2m} \,F = 0 \quad  {\rm in}
\,\,\, \ren,
 }
\end{equation}
so there occurs the linear operator $\BB$ given in \ef{3}. In
addition, the kernel $F$ is defined in such a way that,  in the
sense of distributions (or other suitable $L^p$-type metrics):
 \beq
 \label{DD11}
  \mbox{$
 b(t)*u_0 \to u_0, \,\,\,
 t \to 0^+,
 \quad \mbox{or}
 \quad \int\limits_{\ren} F(y) \chi(t^{\frac 1{2m}} y)\, {\mathrm d}y \to \chi(0)
 \quad \forall \chi  \in C_0^\iy(\ren),
 $}
  \eeq
  justifying initial data. These define the unique
  rescaled kernel $F$.


 On the other hand,
 using
the Fourier transform yields the following equivalent
representation of $F$:
\begin{equation}
\label{FT1}
 {\mathcal F}(b(\cdot,t))(\o) = {\mathrm e}^{- \ii |\o|^{2m}t} \LongA
 {\mathcal F}(F(\cdot))(y) = {\mathrm e}^{- \ii |y|^{2m}}.
\end{equation}
 For $m=1$, this gives the ``Gaussian" exponential profile
 \beq
 \label{Ga1}
  \tex{
 F(y)=  \frac 1{(4 \pi \ii)^{N/2}} \, {\mathrm e}^ {\frac{\ii
 |y|^2}4} \quad (m=1).
 }
 \eeq



 It follows that $F(y)$ is highly oscillatory as $y \to \iy$. In
 particular, for \ef{Ga1}, we have
  \beq
  \label{Ga2}
   \tex{
  F(y) = \frac 1{(4 \pi \ii)^{N/2}}\big[ \cos\big( \frac {|y|^2}4\big) +
   \ii \sin\big( \frac {|y|^2}4\big) \big], \quad
   |F| \equiv {\rm const.}, \,\,\,
  F  \not \in L^p(\ren), \,\,\, p \ge 1.
  }
   \eeq

For arbitrary $m \ge 2$,  the asymptotic behaviour of $F(y)$ for
$|y| \gg 1$ is covered by the classic WKBJ asymptotics. Namely,
fixing in \ef{ODEf} two main leading terms for the radial kernel
 $F=F(z)$, for $z=|y| \to \iy$,
 \beq
 \label{mm1}
  \tex{
 - \ii (-1)^m F^{(2m)} + \frac 1{2m}\, F'z+...=0
}
 \eeq
 yields, in the first approximation, the following exponential asymptotic
 behaviour:
  \beq
  \label{mm2}
  \tex{
   F(y) \sim \eee^{a |y|^\a} \LongA \a= \frac {2m}{2m-1} \andA (\a a)^{2m-1} = \frac{(-1)^{m+1}
   \ii}{2m}.
   }
   \eeq
Hence, there exist $2m-1$ different complex solutions $\{a_k\}$
belonging to a circle in ${\mathbb C}$:
 \beq
 \label{mm3}
 \tex{
  |a_k| = z_m= \frac 1 \a \, (2m) ^{- \frac 1{2m-1}} < 1  \forA m
  \ge 1.
  }
   \eeq
 Obviously, we are interested in those roots $a_k$, for which $\RR
 a_k \le 0$. Otherwise these will be exponentially growing oscillatory
 solutions that will be ``too much" non-integrable.
 It is clear that there exists the purely imaginary root with the
 main asymptotic oscillatory behaviour at infinity:
  \beq
  \label{mm4}
   \tex{
    a_0= z_m \, \ii \LongA F(y) \sim \cos (z_m |y|^\a) + \ii \sin
    (z_m |y|^\a) \asA y \to \iy.
     }
     \eeq

On the other hand, the ODE \ef{ODEf} admits solutions with a power
decay:
 \beq
 \label{mm8}
  \tex{
   \frac 1{2m}\, F' z + \frac N{2m} \, F+...=0, \quad z=|y| \gg 1 \LongA \tilde F(y) \sim
   \frac C{|y|^N} \asA y \to \iy.
    }
     \eeq
 Since $\tilde F(y)$ is then are non-oscillatory (for being used as in
 \ef{DD11}) and, in addition, are ``too much"
   non-integrable as $y \to \iy$, such
 asymptotics are not acceptable for the fundamental kernel $F(y)$.
 We then arrive at the following simple, but interesting and, in fact, a key property of $F(y)$,
 which will affect our analysis (especially, in the ``nonlinear" cases):
  \beq
  \label{mm9}
   \mbox{all asymptotic components of the rescaled kernel $F(y)$ as $y \to \iy$ are
   oscillatory.}
    \eeq

\subsection{Convolution: a unitary group}

Thus, the unique weak solution $ u(x,t)$  of the Cauchy problem
\ef{1} for any data $u_0 \in {\mathcal L}'$
 is given by
 the  Poisson-type integral for $t \in \re$:
\begin{equation}
\label{usol1}
 \tex{
 u(x,t) = b(t)*u_0
 \equiv {\mathrm e}^{-\ii  (-\D)^m \,t}u_0= t^{-\frac N{2m}}
\int\limits_{\ren}F \big((x-z)t^{-\frac 1{2m}}\big)u_0(z) \,
{\mathrm d}z,
 }
\end{equation}
where $\{ {\mathrm e}^{-\ii  (-\D)^m \, t}\}_{t \in \re}$ is the
corresponding unitary group. In what follows, for some
convenience, we take $t>0$ only (that
 implies no trouble in defining the flow in \ef{usol1}), so actually we deal with
 the semigroup $\{ {\mathrm e}^{-\ii  (-\D)^m \, t}\}_{t \ge 0}$.



\section{Discrete real spectrum and eigenfunctions of $\BB$}
 \label{S3}

This section is devoted to some preliminary analysis of spectral
properties of the key pair of linear rescaled operators
$\{\BB,\BB^*\}$ that appear after long-time ($t \to +\iy$) and
short-time ($t \to T^-$) respectively rescaling of the LSE--$2m$
(\ref{1}). In fact, this explains in a reasonably  brief manner
several necessary key properties of the pair to be used later on.
However, we must admit that some of the issues will require  a
hard work to justify by classic theory.
 This will take a full separate paper \cite{GK2mClass}.
 Nevertheless, we hope that
listing a full collection of some involved spectral properties
will be convenient for at least some of the Readers, who are
interested in general understanding of how this approach works and
who do not require full mathematical details.

\subsection{First step to the domain of $\BB$ in a weighted $L^2$-space}

 We now study spectral properties of the first appeared linear operator $\BB$ given in
 \ef{ODEf}
in the space $L^2_{\rho}(\ren)$
 with the exponential weight:
 \begin{equation}
 \label{rho11}
 \tex{
 \rho(y) = {\mathrm e}^{- |y|^\a}>0 \quad {\rm in} \,\,\, \ren
 \whereA \a= \frac{2m}{2m-1}.
 }
 \end{equation}
 By
$\langle \cdot, \cdot \rangle$, we  denote the
 standard $L^2$-product:
  \beq
  \label{sc1}
   \tex{
   \langle v, w \rangle = \int\limits_{\ren} v(y) \overline{w(y)} \,
   {\mathrm d}y.
    }
    \eeq



As customary,  $ H^{2m}_\rho(\ren)$ denotes
  a Hilbert space
 of functions
with the inner product
 \beq
 \label{vsc}
  \tex{
 \langle v,w \rangle_\rho = \int\limits_{\ren} \rho(y) \sum\limits_{k=0}^{2m} D^{k}
 v(y) \, \overline {D^{k} w(y)} \,{\mathrm d} y,
 }
 \eeq
 where $D^k v$ stands for the vector $\{D^\beta v, \,\, |\beta|=k\}$,
 and the norm
 \begin{equation}
 \label{vnorm}
  \tex{
 \|v\|^2_\rho = \int\limits_{\ren} \rho(y) \sum\limits_{k=0}^{2m} |D^{k}
 v(y)|^2 \, {\mathrm d} y.
  }
 \end{equation}
 Obviously,  $
  H^{2m}_\rho(\ren) \subset L^2_\rho(\ren) \supset L^2(\ren)$.
Introducing the weighted Sobolev space  $H^{2m}_\rho(\ren)$ is  a
first step to better understanding
 the necessary and natural domain of $\BB$, as stated in
the proposition below. However,  the space $\LL$ with the
exponentially decaying weight \ef{rho11} is evidently too wide, so
we cannot expect any good spectral properties therein.
Nevertheless, we now prove the following:

\begin{proposition}
 \label{dom11}
 ${\bf B}$ is a bounded linear operator from
 $ H^{2m}_\rho(\ren)$ to  $ L^2_\rho(\ren)$.
 \end{proposition}

 \noi {\em Proof.}
    It follows from
 (\ref{ODEf}) that ${\bf B}v \in L^2_\rho(\ren)$,
    if
 \begin{equation}
 \label{2.3Eg}
  \tex{
 \int\limits_{\ren} \rho(y)|y \cdot \nabla v|^2 \, {\mathrm d}y  \le C \|v\|_\rho^2
 \quad \mbox{for any} \,\,\, v \in  H^{2m}_\rho(\ren),
  }
 \end{equation}
  where $C>0$ is a constant.
Let $h \in C^\infty(\ren)$ be a function such that $h(y)=1$ for $|y| \ge 2$ and $h(y) = 0$
for $|y| \le 1$. Since the inequality
 $$
  \tex{
\int\limits_{\ren} \rho(y)|y \cdot \nabla [(1-h(y)) v]|^2 \,
{\mathrm d}y \le C_1 \|v\|_\rho^2
 }
  $$
  is obvious, it suffices to show that
 $$
  \tex{
\int\limits_{\ren} \rho(y)|y \cdot \nabla (hv)|^2 \, {\mathrm d}y
\le C \|h v\|_\rho^2, }
 $$
 i.e., proving (\ref{2.3Eg}), we can suppose
that $v \in C_0^\infty(\ren)$ vanishes for all $|y| \le 1$.

Let $(r,\theta_1,...,\theta_{N-1})$ be the spherical coordinates in $\ren$. Since
 $
 |y \cdot \nabla v| \le r |v_r|,
 $
 it suffices to verify that
 \begin{equation}
\label{2.4Eg}
 \tex{
 \int\limits_0^\infty r^{N+1} |w_1(r)|^2 {\mathrm e}^{-r^\a} \,
  {\mathrm d}r \le C_2 \int\limits_0^\infty r^{N-1}
 |w_1^{(2m-1)}(r)|^2 {\mathrm e}^{-r^\a}\, {\mathrm d}r,
 }
 \end{equation}
if the left-hand side in bounded,
and apply this estimate with $w_1=v_r$.

Let $q=N-1$ or $q=N-3$, $\g = \a-1 = \frac 1{2m-1}$ and $ k \in
\{0,1,...,2m-2\}$. Then using the inequality
 $$
  \tex{
 \int\limits_0^\infty r^{q+2k \g} |w'(r) + r^\g w(r)|^2 {\mathrm e}^{-r^\a} \, {\mathrm d}r \ge 0,
  }
 $$
integrating by parts again implying that the right-hand side
converges, we obtain  that
 $$
  \tex{
 \int\limits_0^\infty r^{q+2k \g} |w'(r)|^2  {\mathrm e}^{-r^\a} \, {\mathrm d}r \ge \g
  \int\limits_0^\infty r^{q+2k \g +2 \g} |w(r)|^2  {\mathrm e}^{-r^\a} \, {\mathrm d}r.
   }
 $$
Simple iteration implies (\ref{2.4Eg}) with $C_2 = \g^{1-2m}$,
 completing the proof.  $\qed$

 \ssk

The result also follows from a general estimate  in
\cite[Lemma~2.1]{Hei}, which goes back to the Hardy classical
inequality \cite{Hardy}. In a similar (or obvious in (ii)) manner,
introducing the ``adjoint" spaces with the reciprocal weight
 $L^2_{\rho^*}(\ren)$, with the weight \ef{ww*}.
  \beq
  \label{ww*}
   \tex{
    \rho^*(y)= \frac 1 {\rho(y)}= {\mathrm e}^{ |y|^ \a},
    }
    \eeq
 we have the following:

\begin{corollary}
 \label{Col.1}
 $\BB$ is  bounded as an operator

 \noi {\rm (i)} $\BB: \,\,H_{\rho^*}^{2m}(\ren) \to L^2_{\rho^*}(\ren)$ and

\noi {\rm (ii)} $\BB: \,\,H_{\rho^*}^{2m}(\ren) \to
 L^2_\rho(\ren)$.

 \end{corollary}

\ssk

\noi{\bf  Remark for $m=1$.} As customary,  in the  second-order
case $m=1$, there appear some extra possibilities and
``symmetries". Namely, then ${\bf B}$ admits a formal symmetric
representation
 \beq
 \label{ss1}
 \tex{
 {\bf B} =\ii \frac 1  \kappa \, \nabla \cdot (\kappa \n) + \mbox{$\frac N2$}\, I, \quad
 \mbox{with the ``weight"} \,\,\, \kappa(y) =
 {\mathrm e}^{-\ii \frac{|y|^2}4},
  }
  \eeq
 in the  weighted space $L^2_\kappa$ with the complex weight $\kappa \not =0$ and
  hence with an ``indefinite metric"; cf.   Azizov--Iokhvidov \cite{AI89}, which we will need to
   refer to later on a few times at least.
   We do not know any reasonable application of the complex symmetric form
   in \ef{ss1}.
    For instance, as usual, the symmetry \ef{ss1} implies the
    formal
    orthogonality of eigenfunctions:
     $$
     \tex{
     \langle \psi_\b, \psi_\g \rangle_\kappa \equiv
 \int \kappa \psi_\b \psi_\g
      = 0 \forA \b \not =
     \g,
     }
     $$
   but it is not that easy to  find suitable applications of this
   in view of the indefinite metric involved.
   Anyway,
 by no means, we are going to rely on this kind of a pseudo-symmetric
 representation of the operator for $m=1$, especially, since
  for  $m \ge 2$, this illusive complex symmetry of $\BB$ disappears without a trace.
\subsection{Group with the infinitesimal generator $\BB$}

Before introducing detailed spectral properties of $\BB$, we
present a simple derivation of its group
 for proper weak solutions
  to be heavily
 used in what follows.

Thus, the rescaled solution of \ef{1} defined as
 \begin{equation}
 \label{Resc11}
  \mbox{$
  w(y,\t)=t^{\frac N{2m}}
u\bigl(yt^{\frac 1{2m}},t \bigr), \quad \mbox{where} \,\,\, \t =
\ln t \in \re \quad (t>0),
 $}
\end{equation}
 satisfies the
necessary  rescaled equation
\begin{equation}
\label{weq}
 w_\t = {\bf B}w \quad (\mbox{the operator \,$\BB$\, is  as in (\ref{ODEf})}).
\end{equation}
Then $w(y,\t)$ solves the CP for (\ref{weq}) in $ \ren \times
\re_+$ with
 data at $\t=0$ (i.e., at $t=1$)
 \begin{equation}
\label{win0}
   w_0(y) =  u(y,1) \equiv b(y- \cdot,1) * u_0(\cdot).
\end{equation}
Rescaling  convolution (\ref{usol1}) yields
 the following explicit representation of the group with the infinitesimal generator
 $\BB$:
\begin{equation}
\label{geta1}
 \mbox{$
 w(y,\t) =  {\mathrm e}^{{\bf B}\t} u(y,1) \equiv
 \int\limits_{\ren} F\bigl(y- z  {\mathrm e}^{-\frac 1{2m} \, \t}\bigr)
u_0( z )\,{\mathrm d}  z  \quad \mbox{for} \quad  \t \in \re.
 $}
\end{equation}

Performing another
 rescaling
  \beq
  \label{Resc11N}
   \mbox{$
  w(y,\t)=(1+t)^{\frac N{2m}}u \bigl(y(1+t)^{\frac 1{2m}},t\bigr)
  \whereA  \t=\ln(1+t): \re_+ \to \re_+,
   $}
   \eeq
   we
obtain the solution $w(y,\t)$ of the Cauchy problem for equation
(\ref{weq}) with initial data
 $
 w_0(y) \equiv u_0(y).
 $
  Rescaling
(\ref{usol1}), we deduce a more complicated, but standard (without
the relation (\ref{win0})) representation of the
semigroup for $\t \ge 0$
 \begin{equation}
 \label{semNN}
  \mbox{$
 w(y,\t) = {\mathrm e}^{{\bf B}\t} u_0 \equiv  (1-{\mathrm e}^{-\t})^{-\frac N{2m}}
 \int\limits_{\ren}
 F\bigl((y-  z  {\mathrm e}^{-\frac 1{2m} \, \t}) (1-{\mathrm e}^{-\t})^{-\frac 1{2m}}\bigr) u_0( z ) \,{\mathrm d}  z .
  $}
 \end{equation}
By the H\"older inequality (see e.g.,  \ef{jj2} below), it is easy
to see that
 \beq
 \label{ww11}
 w (\cdot,\t) \in \LL \quad \mbox{for all} \quad \t>0 \quad (u_0 \in \LLL).
 \eeq






\subsection{Spectral decomposition of $\BB$: a first step to  discrete spectrum and
 eigenfunctions via converging expansion of convolution}

We now in a position to make more clear a proper definition of the
necessary operator $\BB$ to be exploited  later on.
 We then confirm the actual existence and  the so-called evolution
completeness/closure of some eigenfunctions for initial data $u_0
\in L^2_{\rho^*}(\ren)$. It is worth to stress out how simple such
derivations are here, unlike a more standard spectral and analytic
continuation techniques applied in \cite{GK2mClass}
to justify necessary  hard properties of resolvent poles and
related issues.

 It can be derived
from the ODE (\ref{ODEf}) for the rescaled kernel $F(|y|)$ and
also from \ef{FT1}
 that the higher-order derivatives of $F$
can be  estimated as follows:
\begin{equation}
\label{ff1}
|D^\beta F(y)| \le C \, (1+|y|)^{(\a-1)|\b|} \quad
\mbox{in \,\,$\ren$}.
\end{equation}
 Actually, according to \ef{mm2}, \ef{mm3}, a sharper estimate
 includes the factor $|\a z_m|^l$, i.e.,
  \beq
  \label{mm6}
   \tex{
   |\a z_m|^l = (2m)^{-\frac l{2m-1}}  \to 0 \asA l \to \iy,
    }
    \eeq
    so this improves convergence of the series to appear later
    on.






  Consider  Taylor's power series of the analytic kernel
  $F(\cdot)$ on compact subsets $y \in \{|y|\le L\}$, with an $L \gg 1$,
\begin{equation}
 \tex{
\label{Tay1}
    F\big(y-z{\mathrm e}^{- \frac \t{2m}}\big) =
\sum\limits_{(\b)}{\mathrm e}^{-\frac{|\b|\t}{2m}} \, \frac
{(-1)^{|\b|}}{\beta!} D^\beta F(y) z^\b \equiv
\sum\limits_{(\b)}{\mathrm e}^{-\frac{|\b|\t}{2m}} \frac
1{\sqrt{\b !}} \, \psi_\b(y)z^\b,
 }
\end{equation}
 where $z^\b \equiv  z_1^{\beta_1}...z_N^{\beta_N}$ and $\psi_\b$
 are in fact
normalized eigenfunctions of ${\bf B}$; see below. This series
converges uniformly on compact subsets in $z \in \ren$. Indeed,
for $|\b|=l\gg 1$, we have the following approximate  estimate of
the expansion coefficients:
 \beq
 \label{GG1}
  \tex{
\Big| \sum\limits_{|\beta|=l}\frac {(-1)^l}{\beta!}\, D^\beta F(y)
z_1^{\beta_1}...z_N^{\beta_N}\Big| \le C \frac 1{l!}\,
(1+|y|)^{(\a-1)l} \, (1+|z|)^{l},}
  \eeq
 where we have used a rough bound by Stirling's formula (as usual, we often omit the
 lower-order multipliers
  $\sim C^l$ in \ef{GG1}):
 $$
 \tex{
 \b! \ge \big[ \big( \frac l N \big)!\big]^N \sim N^l l! \, .
 }
 $$





 Finally, we arrive at the following
 representation of the solution:
\begin{equation}
\label{gsol}
 \tex{
 w(y,\t) = \sum\limits_{(\b)} {\mathrm e}^{-\frac{|\b| \t}{2m}}M_{\beta}(u_0)
\psi_{\beta}(y),
 }
\end{equation}
where $\l_\b = -\frac{|\b|}{2m}$ and $\psi_\b(y)$ are the
eigenvalues and eigenfunctions of ${\bf B}$ and
\begin{equation}
\label{Mmom}
 \tex{
 M_{\beta}(u_0) = \frac 1{\sqrt{\b !}} \int\limits_{\ren}
z_1^{\beta_1}...z_N^{\beta_N}u_0(z)\,{\mathrm d} z }
\end{equation}
are the corresponding moments of the initial datum $w_0$ (recall
the relation (\ref{win0}) between $w_0$ and $u_0$). We will show
that, in terms of the dual inner product $\langle \cdot,\cdot
\rangle$ in $L^2(\ren)$,
\begin{equation}
\label{MBB1} M_\b(u_0)= \langle w_0, \psi_\b^*\rangle,
\end{equation}
where $\{\psi_\b^*\}$ are
 polynomial eigenfunctions of the adjoint
operator ${\bf B}^*$ to be described in greater detail in Section
\ref{S4}. It is not difficult to check that \ef{gsol}
 uniformly converges on any compact subset $y \in \{|y|\le L\}$,
 since for $l=|\b| \gg 1$,
 \beq
 \label{jj1}
 \tex{
 |M_{\beta}(u_0)
\psi_{\beta}(y)| \le C  \frac 1{l!} \,  L^{(\a-1)l}
\int\limits_{\ren} (1+|z|)^l |u_0(z)|\, {\mathrm d}z,
 }
 \eeq
 where we estimate the last integral by using the radial
 variable and the H\"older inequality,
  \beq
  \label{jj2}
    \begin{matrix}
   \int\limits_{\ren}
(1+|z|)^l |u_0|\, {\mathrm d}z = \int\limits_{\ren} (1+|z|)^l
 \frac 1{\sqrt{\rho^*}} \, \sqrt{\rho^*}|u_0|\, {\mathrm d}z
  \ssk\ssk\ssk\\
  \le \, \sqrt{\int_{\ren} (1+|z|)^{2l}
 \frac 1{{\rho^*}}\, {\mathrm d}z}\,\,
 \sqrt{\int_{\ren} {\rho^*}|u_0|^2\, {\mathrm d}z}.
   \end{matrix}
  \eeq
 The last integral is bounded since $u_0 \in L^2_{\rho^*}$ by
the  assumptions, while the first integral for very large $l$ can
be roughly estimates
 by Stirling's formula as follows:
 \beq
  \label{jj3}
   \tex{
   \sqrt{\int_{\ren} (1+|z|)^{2l}
 \frac 1{\rho^*}\, {\mathrm d}z} \sim \sqrt{\int_0^\iy r^{N-1+2l}
 {\mathrm e}^{-r^\a}\, {\mathrm d}r} \sim  \big( \frac {N+2l}{\a
 \, {\mathrm e}
 }\big)^{\frac{N+2l}{2\a}} \sim
 \big( \frac{2}{\a \,{\mathrm e}} \big)^{\frac l \a}
  l^{\frac l \a}.
  }
  \eeq
Since $\a>1$, the right-hand side is essentially  smaller
  than $l! \sim l^l$
 in the denominator in \ef{jj1}. Hence,
this series converges uniformly and to an analytic solution as
expected, and not surprisingly for the LSE.




 Recall that \ef{semNN}
 gave  the actual corresponding group $\{T(\t)={\mathrm e}^{{\bf
B}\t}\}_ {\t \in \re}$. Here, we can apply the same expansion
analysis as above, which directly determines the  eigenfunctions
of ${\bf B}^*$. First, by Taylor's expansion in the $F$-term,  we
obtain
 \begin{equation}
 \label{w1exp}
  \tex{
 w(y,\t) = \sum\limits_{(\mu)} \frac {(-1)^{|\mu|}}{\mu !} D^\mu F\big(y (1-{\mathrm e}^{-\t})^{-\frac1{2m}}
 \big)
 ({\mathrm e}^\t-1)^{-\frac{|\mu|}{2m}} \int\limits_{\ren} z^\mu w_0(z) \, {\mathrm d}z.
  }
 \end{equation}
Second,  applying  Taylor expansions for the $f$-terms in
(\ref{w1exp}),
 \beq
 \label{kam1}
  \tex{
 F\big(y (1-{\mathrm e}^{-\t})^{-\frac1{2m}}\big) = \sum_{(\nu)} \frac 1{\nu !} D^\nu
 F(0) y^\nu \, (1-{\mathrm e}^{-\t})^{-\frac{|\nu|}{2m}}
 }
 \eeq
 and for functions   $(1-{\mathrm e}^{-\t})^{-\frac 1{2m}}$
 and $(1-{\mathrm e}^{-\t})^{-\frac N{2m}}$
 in terms of ${\mathrm e}^{-\frac {k \t}{2m}}$, $k=0,1,...$,
 we arrive at a similar
 representation of the semigroup
\begin{equation}
\label{w1New}
 \tex{
 w(y,\t) = {\mathrm e}^{{\bf B}\t}w_0   \equiv \sum_{(\b)}
{\mathrm e}^{-\frac{|\b|\t}{2m}} \tilde M_\b(w_0)  \psi_{\b}(y),
\quad \t \ge 0,
  }
\end{equation}
where the expansion coefficients $\tilde M_\b(w_0)$ are
  dual products $\langle w_0,  \psi^*_\beta \rangle$ with the
polynomial eigenfunctions $\psi_\beta^*$ of the adjoint operator
${\bf B}^*$. Expansion (\ref{w1New})  determines the adjoint
eigenfunctions $\{\psi_\b^*\}$ that are ``orthogonal" to
$\{\psi_\b\}$ in $L^2(\ren)$ (all products and metrics, with the
convention \ef{psiNN}, to be introduced later on). In Section
\ref{S4} devoted to the adjoint operator ${\bf B}^*$, we perform a
simpler derivation of explicit formulas for polynomials
$\{\psi_\b^*\}$.


\ssk

Thus, using the expansion \ef{w1New}, we are now in a position to
present the first proper definition (via its spectral
decomposition) of our bounded linear operator in
$H^{2m}_{\rho^*}(\ren)$, which we denote by ${\mathbb B}$:
 \beq
 \label{de1}
 \tex{
 {\mathbb B} w_0= \frac{\mathrm d}{{\mathrm d}\t}\, {\mathrm e}^{{\mathbb B}\t}w_0
 \big|_{\t=0}   \equiv \sum_{(\b)}\l_\b
 \tilde M_\b(w_0)  \psi_{\b}(y) \forA w_0 \in H^{2m}_{\rho^*}(\ren).
 }
 \eeq
 Actually, ${\mathbb B}$ can be considered as a restriction of a
 more general operator $\BB$ defines in a so-called  wider space of
 closures. Since this extension is not essential for our main
 applications and ${\mathbb B}$ is sufficient, we postpone this
 rather technical procedure until Appendix D at the paper end.

 \ssk

 \noi {\bf Comment on ``extended eigenfunctions".} With the definition of our operator
  \beq
  \label{BB11}
  {\mathbb B}: H^{2m}_{\rho^*}(\ren) \to \LL
   \eeq
 by its expansion \ef{de1}, we face the following technical
 difficulty. Namely, the eigenfunctions $\{\psi_\b\}$, which
 actually generate such a ${\mathbb B}$, {\em do not belong} to its
 domain in \ef{BB11}, and we refer to them as to {\em extended
 eigenfunctions} (i.e., as we will show in Appendix D, belonging to an extended space).
  However, for simplicity, we continue to call
 them simply {\em eigenfunctions}, bearing in mind that, in
 Appendix D, we are going to construct an extended space of
 closures, where all $\psi_\b$ belong to. We then  restore the
 original operator $\BB$ rather than its restriction ${\mathbb
 B}$, though, as we  mentioned, for our main PDE
 applications, present ``spectral theory" of the restriction
 ${\mathbb B}$ is more than sufficient.

\ssk

We now summarize the above results concerning the introduced
operator ${\mathbb B}$. Recall that, originally, ${\mathbb B}$ is
defined by the rescaled convolution
 \ef{semNN}, with the corresponding space and
the domain. For further applications, we restricted ${\mathbb B}$
and defined it
 via the spectral decomposition such as \ef{de1}, which  demanded
  special topologies for a proper convergence.


\begin{proposition}
 \label{Pr.BB}

   The semigroup expansion series
 $(\ref{w1New})$ for $w_0 \in L^2_{\rho^*}(\ren)$, and
  the eigenfunction expansion of ${\mathbb B}$  in
 $(\ref{de1})$ for $w_0 \in H^{2m}_{\rho^*}(\ren)$
 converge:

 \ssk

 {\rm (i)} uniformly on compact subsets in $y$, and

\ssk

 {\rm (ii)} in the
 mean in $\LL$.
 \end{proposition}

\noi{\em Proof.} (i) has been proved. (ii) The convergence in the
mean in $\LL$ is not straightforward for such a bad oscillatory
and growing basis functions in $\Phi $. Therefore, estimating the
terms in \ef{gsol}, we need to include both decaying multipliers
in the right-hand side of \ef{jj3}, so that the following estimate
of a typical integral is key by using Stirling's asymptotics of
the Gamma function:
 \beq
 \label{ffg1}
  \begin{matrix}
 l^{2l(\frac 1 \a-1)} \big( \frac 2{\a \, {\mathrm
 e}}\big)^{\frac{2l}\a} \int\limits_0^\iy z^{N-1} {\mathrm e}^{-z^\a}
 z^{2(\a-1)l}\, {\mathrm d}z \sim
 l^{2l(\frac 1 \a-1)} \big( \frac 2{\a \, {\mathrm
 e}}\big)^{\frac{2l}\a} \int\limits_0^\iy s^{\frac {N+2(\a-1)l}\a-1} {\mathrm e}^{-s}
 \, {\mathrm d}s \quad\qquad \ssk\ssk\\
 \sim \, l^{2l(\frac 1 \a-1)} \big( \frac 2{\a \, {\mathrm
 e}}\big)^{\frac{2l}\a}
  l^{\frac {2( \a-1)l}\a} \big( \frac {2(\a-1)}{\a \, {\mathrm
 e}}\big)^{\frac{2(\a-1)l}\a}
 = (\a-1)^{\frac{2(\a-1)l}\a} \big( \frac 2{\a \, {\mathrm
 e}}\big)^{{2l}} < \big( \frac 4{{\mathrm e}^2}\big)^l. \qquad
   \end{matrix}
   \eeq
 Since $\a \in (1,2]$, this exponentially (but not superexponentially as before for
 the coefficients of \ef{gsol} on compact subsets in $y$) decaying
 coefficients guarantee the convergence of the series in $\LL$.
 The factor in \ef{mm6} is then not necessary for
 convergence. $\qed$





\subsection{Discrete spectrum}

Thus, the series \ef{de1}
 is a {\em spectral decomposition} (actually, the eigenfunction expansion
 for a discrete spectrum) of our linear non self-adjoint
operator, which we have denoted by ${\mathbb B}$.
 Hence, we  treat it as a linear bounded operator given in
  \ef{BB11},
   when the eigenfunction expansion \ef{w1New} is defined in two
   topologies: local uniform convergence, and strong convergence
   in the mean in $\LL$.
   In this case, the expansion in $\LLL$,
   with the same convention on topologies (see \ef{psiNN} explaining why
   $\bar \psi_\b^*$ occur in the last products),
    \beq
    \label{BB12}
     \tex{
    w \in \LLL \LongA w= \sum c_\b \psi_\b, \quad c_\b= \langle w,
    \bar \psi_\b^* \rangle,
    }
     \eeq
is naturally treated as the eigenfunction series representation of
the embedding operator $\LLL \to \LL$ in terms of the
eigenfunctions $\Phi $ of ${\mathbb B}$.

For convenience and more systematic understanding the new class of
linear operators introduced,
 we
clearly state necessary spectral properties of the operator
${\mathbb B}$. According to \cite{GK2mClass},
 the present operator ${\mathbb B}$,
which was not still properly defined and has been
uniquely characterized by its ``spectral eigenfunction
decomposition" \ef{de1} only, is not allowed to have the positive
part of the point spectrum with polynomial eigenfunctions, as
shown in \cite{GK2mClass},
 where this
is done by a direct introducing ``radiation conditions" posed at
infinity.

 Thus, according to our operator definition \ef{de1}, we ascribe
 to ${\mathbb B}$ a countable set of numbers, which, for
 convenience, we continue to call its discrete spectrum:
 \begin{equation}
\label{spec1}
 \fbox{$
 \tex{
 \sigma({\mathbb  B})=
\big\{\lambda_\b = -\frac{|\b|}{2m}, \,\, |\b| = 0,1,2,...\big\},
 }
  $}
\end{equation}
where  eigenvalues $\l_\b$ have finite multiplicity with
eigenfunctions\footnote{Actually, {\em extended eigenfunctions},
since $\psi_\b \not \in H^{2m}_{\rho^*}(\ren)$; a standard meaning
{\em eigenfunctions} of $\BB$ to be restored in Appendix D by
introducing  an {\em extended space of closures of finite
eigenfunction expansions}.}
\begin{equation}
\label{eigen}
 \fbox{$
 \tex{
\psi_\beta(y) = \frac{(-1)^{|\beta|}}  {\sqrt {\beta !}}
  D^\beta F(y) \equiv \frac{(-1)^{|\beta|}}{\sqrt{\beta !}}
\big(\frac {\partial}{\partial y_1}\big)^{\beta_1}... \big(\frac
{\partial}{\partial y_N}\big)^{\beta_N} F(y).
 }
  $}
\end{equation}
The existence of such eigenvalues and eigenfunctions
 is dictated by \ef{de1}. It is worth mentioning that the
 same
 follows by applying
 $D^\beta$ to the elliptic equation (\ref{ODEf})
 (here $\BB$ stands for its differential expression, so we are not obliged to use ${\mathbb
 B}$): for any $\b$,
\begin{equation}
\label{idB}
 \tex{
 D^\beta{\mathbf B}F \equiv {\mathbf B}D^\beta
F+  \mbox{$\frac{|\b|}{2m}$} \, D^\beta F=0
 \LongA \big\{D^\b F,\l_\b=- \frac {|\b|}{2m}\big\}
 \,\,\mbox{is a
  pair for $\BB$}.
  }
\end{equation}

Let us fix some other properties of such extended eigenfunctions:

\begin{lemma}
\label{lemspec}


{\rm (i)} The subset of eigenfunctions $\Phi =\{\psi_\b\}$ is
complete
  in $L^2_{\rho}(\ren)$, and

 {\rm (ii)}  $\Phi $ is $L^2_{\rho^*}$-evolutionary closed
 in the sense that the eigenfunction expansion $(\ref{w1New})$,
 which converges in the means and uniformly on compact subsets,
presents  the rescaled solution of the LSE $(\ref{1})$ for any
data $u_0 \in L^2_{\rho^*}(\ren)$; and


\end{lemma}

\noi{\bf Remark: towards more general integral evolution equations
and rescaled operators ${\mathbb B}$.} As we have mentioned, all
the above results can be justified by classic spectral methods
associated with the given rescaled differential operators,
\cite{GK2mClass}.

 However, it is worth mentioning now that all the (i)--(iii)
remain valid for more general class of {\em integral evolution}
(pseudo-differential)  equations \ef{geta1}, where
 \beq
 \label{gener11}
  \fbox{$
 F(y) \,\,\, \mbox{is an arbitrary sufficiently good analytic kernel.}
 $}
  \eeq
In other words, \ef{geta1} then do not correspond to any linear
PDE.
 Furthermore, the expansions
 \ef{Tay1} and \ef{gsol} can be also prescribed. Overall, this
 gives the spectral results similar to those in (i)--(iv), where
 extra efforts to justify the functional topology required are
 necessary.

In this connection, it is key to emphasize that, in the present
most general situation, no extra powerful tools of spectral
theory, developed in \cite{GK2mClass}
 for the present
Schr\"odinger operators, will be at hand. Then,  we will be
inevitably attached to a different functional framework, with no
visually available  operator ${\mathbb B}$ (and the ``adjoint" one
${\mathbb B}^*$).

\ssk


Returning to Lemma \ref{lemspec}, note that, as an important
characterization of the eigenfunction set $\Phi $, due to
\ef{mm9}, {\em all of the eigenfunctions}
 $\psi_\b(y)$ satisfy the property \ef{mm9}, i.e., these do not
 have any ``non-fast-oscillatory" asymptotic component as $y \to \iy$
 (since \ef{mm8} has been excluded from the kernel $F(y)$).
Observe also  from \ef{eigen} that, with a proper definition of
such {\em extended (generalized) linear functionals}, to be done
in Appendix D,
 the following can be interpreted as being correct:
\begin{equation}
 \tex{
\label{f1n}
 \langle \psi_0, \psi_0^* \rangle_*=
 \int\limits_{\ren} \psi_0(y) \, {\mathrm d}y =
\int\limits_{\ren}F(y) \, {\mathrm d}y = 1, \quad \mbox{since}
\quad  \psi_0 = F, \,\,\, \psi_0^*=1.
 }
\end{equation}
Recall that, in the usual sense, such oscillatory integrals are
not properly defined. On the other hand, it also follows from
(\ref{eigen}) that, again in a proper extended linear functional
sense (see Appendix D),
\begin{equation}
 \tex{
\label{f2n}
  \langle \psi_\b, \psi_0^* \rangle_*=
 \int\limits_{\ren} \psi_\b(y) \, {\mathrm d}y  = 0
\quad \mbox{for any} \,\,\, |\b| \ge 1.}
\end{equation}
This is easier to believe in view of the  integration by parts,
though the oscillatory integrals are not well defined as well.  In
fact, these equalities express the orthogonality of any $\psi_\b$
to the first adjoint eigenfunction $\psi_0^*=1$ via the dual inner
product, defined as extensions of
 linear functionals prescribed in \ef{BB12}; see again Appendix D.
   The adjoint
eigenfunctions are polynomials which form a complete subset in
$L^2_{\rho}(\ren)$ with the same decaying exponential weight
 \ef{rho11};
 see Section \ref{S4}.

In the second-order case $m=1$, using the rescaled kernel \ef{Ga1}
in \ef{eigen} gives the corresponding {\em generalized Hermite
polynomials} $H_\b(y)$ (given up to normalization constants) via
the generating formula:
  \beq
  \label{kam2}
   \tex{
 \psi_{\beta}(y)=  \frac 1 {(4\pi \ii)^{N/2}} D^\b {\mathrm e}^{\ii \frac{|y|^2}4}
  \equiv   H_{\beta}(y) {\mathrm e}^{\ii \frac{|y|^2}4}.
 }
 \eeq
Note that, in this case, the generalized Hermite polynomials are
obtained from the classic ones \cite[p.~48]{BS} by the change
 $$
 H_\b(y)=H_\b^{\rm class.}(\ii y) \quad(m=1).
 $$



  \smallskip

\noindent {\em Proof of Lemma $\ref{lemspec}$.}

\noindent{\rm (i)}
 \underline{\em Completeness in $L^2_{\rho}$-space}.
    In order to
prove completeness in the metric of $L^2_{\rho}$, as in
\cite[\S~2]{Eg4},
we suppose that there exists  some function $G$ (say,  $G \in
L^2$), which
is orthogonal relative to the inner product in $L^2_{\rho}$ to all
eigenfunctions, i.e.,
 $$
 \mbox{$
  \int \rho(y) D^\alpha F(y) G(y)  \,{\mathrm d}y =0 \quad \mbox{for all} \,\,\,
 \alpha.
  $}
 $$
Since $F$ is analytic, it implies that
 $$
  \mbox{$
  \int \rho(y) F(y-x)  G(y)  \, {\mathrm d}y =0 \quad \mbox{for all $x\in
\ren $.}
   $}
   $$
     Consider the Cauchy problem for
 (\ref{1}) with initial data
 $$
 u_0(x)= \rho(x) G(x) \quad \mbox{in} \quad \ren \quad \big(\mbox{note that} \,\,\,
 \rho G \in L^2_{\rho^*}(\ren)\big).
  $$
One can see from the Poisson-type integral (\ref{usol1})
that
 the solution exists  for all $t \in (0,1]$.
 Then $u(x,t)$ is analytic in $x$.
We have
 $$
  \mbox{$ u(x,1)= \int F(x-y) G(y) \rho(y) \,{\mathrm d}y.
   $}
  $$
   Therefore,
$u(x,1) \equiv 0$. It follows from the standard uniqueness theorem
(see \cite{Hao07} as a guide)
that $u(x,0)=0$, and $G=0$.

\smallskip


\noindent (ii) \underline{\em Evolution closure} in $L^2_{\rho^*}$
has been already proved while studying the convergence of the
series \ef{gsol} and \ef{w1New}.



 This completes the
proof of Lemma \ref{dom11}. $\qed$



\ssk

\noi{\bf Remark: positive part of the ``point spectrum".}
Formally, in addition to \ef{spec1}, there exists another positive
part of the ``real spectrum" of the differential form ${\bf B}$
\ef{ODEf}:
 \beq
 \label{pos1}
  \begin{matrix}
  \s_+({\mathbb B})=
\big\{\lambda_\b = \frac{N+|\b|}{2m}, \,\, |\b| = 0,1,2,...\big\}
\whereA \ssk\ssk\ssk\ssk\\
   \psi_\b^+(y)= \frac 1{\sqrt \b!}\, y^\b+...
\,\,\,\mbox{are $|\b|$th-degree polynomials}.
 \end{matrix}
\end{equation}
We will show how to construct such polynomial eigenfunctions in
Section \ref{S4}, where these are actual eigenfunctions of
 the adjoint operator ${\mathbb B}^*$.
  Though such
``pseudo-Hermite" polynomials $\Phi_+= \{\psi_\b^+ \in
L^2_{\rho}\}$ actually exist, they are not  eigenfunctions since
are not available in the given expansions \ef{w1New}
 as the definition of the operator. In other words, these do not
 satisfy the necessary {\em radiation conditions} at infinity,
 which turn out to be rather non-standard and unusual,
 \cite{GK2mClass}.


\section{\bf Spectrum and  polynomial eigenfunctions of the adjoint operator ${\mathbb B}^*$}
\label{S4}

\subsection{Indefinite metric and domain of the bounded operator ${\mathbb B}^*$}

Using the results obtained above  for differential expression
$\BB$ and its proper restriction ${\mathbb B}$ in \ef{de1}, we now
describe in detail the eigenfunctions of a restriction of the
``adjoint" operator \ef{RO1}, which is  still a differential
expression and we use the notation $\BB$ and not ${\mathbb B}^*$.
 $\BB^*$ will be obtained via blow-up rescaling \ef{RR1}.
 Note that $\BB^*$  is not adjoint in the standard dual metric \ef{sc1}  of
 $L^2$, since, as we have seen,
  \beq
  \label{hh2}
   \tex{
  \hat \BB^*=(\BB)^*_{L^2}=- \BB + \frac N{2m} \, I \equiv  {\bar \BB}^*.
 }
   \eeq
Curiously, $\BB^*$ is then the standard adjoint to $\BB$ (as a
differential form, in $C^\iy_0(\ren)$) in the indefinite complex
metric (cf. \ef{sc1})
 \beq
  \label{sc2}
   \tex{
   \langle v, w \rangle_* = \int\limits\limits_{\ren} v(y){w(y)} \,
   {\mathrm d}y \LongA (\BB)^*_*= \BB^*,
    }
    \eeq
    where the complex conjugation in the second multiplier is not
    assumed.

In fact, being treated in the metric of $L^2$, without a weight as
it used to be, such an indefinite metric is not that challenging
and even can be partially get rid of for the practical use of
eigenfunction expansions; see below. However, some comments are
necessary. Firstly, the set of real $L^2$-functions
 $$
 E_+=\{v \in L^2: \quad {\rm Im}\, v = 0\}
  $$
  is a {\em positive  lineal} (a linear manifold in the field of real numbers)) of the metric,
  i.e.,
\begin{equation*}
\langle v,v \rangle_* > 0 \quad \mbox{for} \quad v \in E_+ \subset
L^2,
  \quad v \not = 0.
\end{equation*}
 The purely imaginary functions,
  $$
  E_-=\{v \in L^2: \quad {\rm Re} \, v =0\},
 $$
 define the corresponding {\em negative lineal}.
Therefore, $L^2
 $ with this metric is {\em
decomposable}:
\begin{equation*}
\mbox{$ v=v_+ + v_- \equiv \frac{v(y)+\bar v(y)}{2} +
\frac{v(y)-\bar v(y)}{2}, \quad \text{where} \quad v_\pm \in E_\pm
\, \Longrightarrow \, L^2 = E_+ \oplus E_-. $}
\end{equation*}
This defines the  corresponding positive {\em majorizing} metric
as follows:
\begin{equation*}
|\langle v,v\rangle_*| \le [v,v]_* \equiv  \langle
v_+,v_+\rangle_\ast -  \langle v_-,v_-\rangle_*,
\end{equation*}
etc. It should be noted that such a  case of the decomposable
space with an indefinite metric having a simple majorizing one is
treated as rather straightforward; see Azizov--Iokhvidov
\cite{AI89} for linear operators theory in spaces with indefinite
metrics\footnote{Basic results of linear operator theory in spaces
with indefinite metrics can be found in Azizov--Iokhvidov's
monograph \cite{AI89}. It was in 1944, when L.S. Pontryagin
published the pioneering  paper ``Hermitian operators in spaces
with indefinite metric" \cite{Pont44}. A new area of operator
theory had been formed from Pontryagin's studies, which, during
the time of the WWII,
 were originated and associated  with some missile-type military research
  (a comment by Yu.S.~Ledyaev).
 This work set by
Pontryagin was continued from 1948 and in the 1950s by M.G.~Krein
\cite{Krein48, Krein50}, I.S.~Iokhvidov \cite{Iokh49}, and
others.}.


    Though, as we have mentioned, linear operator theory in spaces with
    indefinite metrics
    exists
    for more than half a century, we do not think that the complex
    indefinite metric in \ef{sc2} creating the necessary pair
    $\{\BB,\BB^*\}$ (the operator and its adjoint) can play any
    role in what follows. On the other hand, as customary, using
    the metric \ef{sc2} is not that suspicious, since it is
    necessary {\em only} for calculating the expansion
    coefficients according to the standard rule:
     \beq
     \label{sc22}
      \tex{
     v= \sum c_\b \psi_\b \LongA c_\b= \langle v, \psi_\b^*
     \rangle_*,
     }
     \eeq
     while all convergence calculus can be continued to be
     performed in standard metrics.
However, to avoid possible future accusations of using
non-approved indefinite metrics, we are now back to standard
scalar products by noting the following.
 Since $\BB^*$ is shown to have a real point spectrum only, using
the standard $L^2$-metric instead of \ef{sc2} will only mean
replacing the eigenfunctions as follows:
 \beq
 \label{psiNN}
 \psi_\b^* \mapsto \bar \psi_\b^*,
  \eeq
  and we are assuming using  this convention
  any time when necessary and convenient.

In the second-order case $m=1$, \ef{RO1} has a formal complex
symmetric representation
\begin{equation}
 \tex{
\label{Bsymm} {\bf B}^* = \ii  \frac 1{\kappa^*}\, \nabla \cdot
(\kappa^* \nabla ), \quad \kappa^*(y) = {\mathrm e}^{\ii
\frac{|y|^2}4},
 }
\end{equation}
 though we do not use this. Similar to $\BB$, we do not know any advantages, which
  this symmetry in such an indefinite metric can provide.
 However,  for  $m \ge 2$, any
 formal additional symmetry is not available.

For any $m \ge 1$,
 we again consider ${\bf B}^*$ in the weighted
space $L^2_{\rho}(\ren)$ with the same exponentially decaying
weight \ef{rho11},
 and ascribe to ${\bf B}^*$ the domain
 $H^{2m}_{\rho}(\ren)$, which is dense in $L^2_{\rho}(\ren)$.
 As in the previous section,
 $
 {\bf B}^*: \,\, H^{2m}_{\rho}(\ren) \to L^2_{\rho}(\ren)
 $
 is shown to be a bounded linear operator.



\subsection{Semigroup with infinitesimal generator $\BB^*$}

In order to construct the semigroup with the infinitesimal
generator ${\bf B}^*$, we use the rescaled variables corresponding
to blow-up as $t \to 1^-$,
  \beq
  \label{RR1}
   \tex{
  u(x,t) =
w(y,\t), \quad y=\frac x{(1-t)^{1/{2m}}}, \quad \t =
-\ln(1-t):(0,1) \to \re_+.
 }
  \eeq
  Then $w$ solves the problem
\begin{equation}
\label{W1} w_\t= {\bf B}^* w \quad {\rm for} \,\,\, \t >0, \quad
w(0)= u_0 \in L^2_{\rho^*}(\ren).
\end{equation}
Rescaling  (\ref{usol1}), we obtain the following explicit
representation of the semigroup:
\begin{equation}
\label{w12}
 \mbox{$
 w(y,\t) = {\mathrm e}^{{\bf B}^*\t} u_0 \equiv (1-{\mathrm e}^{-\t})^{-\frac
 N{2m}}\int\limits_{\ren}
  F\bigl((y {\mathrm e}^{-\frac 1{2m} \, \t}- z )(1-{\mathrm e}^{-\t})^{-\frac 1{2m}}\bigr)
   u_0( z ) \,{\mathrm d}  z .
   $}
\end{equation}


 \subsection{Spectral decomposition and definition of $\BB^*$:
  using explicit representation of the semigroup}

Similar to $\BB$ in Section \ref{S3}, the original rescaled ``adjoint" operator $\BB^*$
is defined by the convolution \ef{w12}. For the purpose of applications, we will
need its restriction defined in terms of its spectral decomposition
obtained via
the semigroup representation (\ref{w12}).
Comparing semigroups (\ref{w12}) and (\ref{semNN}), we see that
the only difference is in the argument of the rescaled kernel
$F(\cdot)$. Therefore, instead of (\ref{kam1}), we have to use the
following expansion:
  \beq
  \label{Tay4NN}
   \begin{matrix}
  F \bigl((y \, {\mathrm e}^{-\frac 1{2m} \, \t}- z )
  (1-{\mathrm e}^{-\t})^{-\frac 1{2m}}\bigr)=
 \sum_{(\g)} \frac {D^\g F(0)}{\g !}
 \bigl(y {\mathrm e}^{-\frac 1{2m} \, \t}- z \bigr)^\g \, (1-{\mathrm e}^{-\t})^{-\frac {|\g|}{2m}}
, \qquad
 \ssk\ssk\\
\mbox{where} \quad  \mbox{$
 \bigl(y \, {\mathrm e}^{-\frac 1{2m} \, \t}- z \bigr)^\g =
  \sum_{(0 \le \d \le \g)}\, C_\g^\d \,
{\mathrm e}^{-\frac {|\g-\d|}{2m} \, \t}\,
   y^{\g-\d}\, (-z)^\d.
    $} \qquad
     \end{matrix}
  \eeq
 Then, using both (\ref{Tay4NN})
     in (\ref{w12})
     yields
  \beq
   \label{Tay6NN}
 \begin{matrix}
{\mathrm e}^{\BB^* \t} u_0= \sum\limits_{(s \ge 0)} \sum\limits
_{(\g)}
   \sum\limits_{(0 \le \d \le \g)}
{\mathrm e}^{-(\frac {|\g-\d|}{2m}+s) \, \t} \qquad \qquad
\smallskip
\smallskip\smallskip\\
  \times
    \, (-1)^{|\g|}\, y^{\g-\d} \kappa_s\bigl(\frac {|\g|+N}{2m}\bigr)
\frac {1}{(\g-\d) !} \, D_\zeta^{\g-\d} \bigr[\int \bigl(\frac
{1}{\d !} \, D^\d F(\zeta) z^\d\bigr) \, u_0(z)\, {\mathrm
d}z\bigl] \big|_{\zeta=0}.
 \qquad \qquad
 \end{matrix}
 \eeq

This is the expansion over the point spectrum of $\BB^*$,
 \beq
 \label{PP1}
  \mbox{$
 w(y,\t) = {\mathrm e}^{\BB^* \, \t} w_0= \sum\limits_{(\b)} {\mathrm e}^{- \frac{|\b|}{2m} \, \t}
 M_\b^*(u_0) \psi^*_\b(y),
  $}
  \eeq
  where $\psi_\b^*(y)$ are finite polynomial eigenfunctions (see their direct derivations below)
   and
  the expansion coefficients are
  $$
   \mbox{$
  M_\b^*(u_0)= \langle  u_0, \psi_\b \rangle.
   $}
   $$
Convergence of the series \ef{PP1} is studied as in Proposition
\ref{Pr.BB}. Similar to \ef{de1}, the group representation
\ef{PP1} defines the necessary operator ${\mathbb B}^*$ satisfying
\ef{BB11} as
 \beq
 \label{de2}
 \tex{
 {\mathbb B}^* w_0= \frac{\mathrm d}{{\mathrm d}\t}\, {\mathrm e}^{{\mathbb B}^*\t}w_0
 \big|_{\t=0}   \equiv \sum\limits_{(\b)} \l_\b \,
 \tilde M_\b^*(w_0)  \psi_{\b}^*(y), \quad w_0 \in H^{2m}_{\rho^*}(\ren).
 }
 \eeq

Similarly to the case of the operator ${\mathbb B}$
 at the beginning of the proof of Lemma \ref{lemspec},
 in view of standard
regularity properties of linear parabolic flows such as
(\ref{W1}), the semigroup expansion (\ref{PP1}) reveals some key
auxiliary spectral properties of ${\mathbb B}^*$:

(i) the point spectrum is $\s({\mathbb B}^*)= \{ \l_\b =-
\frac{|\b|}{2m} \}$, with any $\l_\b$ having finite multiplicity;

(ii) by the definition, there is no continuous spectrum; and

(iii) polynomial eigenfunctions\footnote{Again, there are {\em
extended} ones, $\not \in H^{2m}_{\rho^*}$, and will restore a
usual meaning in Appendix D.} $\{\psi^*_\b(y)\}$ are closed in
$L^2_{\rho}$, etc.

\ssk

In addition, we have to observe that, unlike \ef{mm9} for
${\mathbb B}$, for the adjoint operator ${\mathbb B}^*$, the
opposite characterization of all the eigenfunctions is in use: all
the eigenfunctions
 \beq
 \label{mm10}
  \mbox{$\psi_\b^*(y)$ of ${\mathbb B}^*$
 are not oscillatory and are of a ``minimal" growth as $y \to \iy$ .}
  \eeq
The last issue of $\psi_\b^*(y)$ being of a ``minimal" growth
(since there are many faster other asymptotics that are
oscillatory) will be key in the nonlinear setting for the QLSE
\ef{QQ1}. In the linear case $n=0$, all those notions admit a
natural (but still not that easy) standard treatment, so we do not
need to stress upon such issues in what follows.

\ssk

   Thus, the definition \ef{de2} of the operator ${\mathbb B}$
   justifies the necessary part of the point spectrum with $\l_\b=
   - \frac{|\b|}{2m}$ and excludes its ``positive part"
   (obviously nonexistent in \ef{de2})
  \beq
   \label{posSp1}
 \tex{
 \l_\b^+=  \frac {N+|\b|}{2m}, \quad |\b| \ge 0,
 \quad \mbox{with}
 \quad \psi_\b^{*+}= \bar \psi_\b.
 }
 \eeq
Nevertheless, for convenience, we will continue to refer to
\ef{mm10} as a simple, efficient, and actually true way for a
correct characterization of necessary eigenfunctions.
 In other words,
 the definition
 of ${\mathbb B}^*$ actually includes special ``radiation-like" conditions at
 infinity, which delete the non-desirable positive spectrum.

\subsection{Discrete spectrum and Hermitian polynomial eigenfunctions}

Thus, defining ${\mathbb B}^*$ by \ef{de2}, with the discrete
spectrum only:
\begin{equation}
\label{SpecN}
 \tex{
  \s({\mathbb B}^*) = \s({\mathbb B})= \big\{ \l_\b = -\frac {|\b|}{2m},
\,\,\, |\b|=0,1,2,...\big\},
 }
\end{equation}
where all eigenvalues have  finite multiplicity, and polynomial
extended eigenfunctions to be determined explicitly shortly. Next,
 as customary,  we
fix  other properties of the adjoint operator in a manner
 similar to  Lemma \ref{lemspec}.

\begin{lemma}
\label{lemSpec2}
 Under the above hypothesis and conditions:

 {\rm (i)} (Extended) eigenfunctions $\psi^*_\b(y)$ are polynomials of
order $|\b| \ge 0 $;

{\rm (ii)} The subset of eigenfunctions $\Phi^*=\{\psi^*_\beta\}$
is complete and closed in $L^2_{\rho}(\ren)$; and

 {\rm (iii)}  $\Phi^*$ is $L^2_{\rho^*}$-evolutionary closed
 in the sense that the eigenfunction expansion $(\ref{PP1})$,
 which converges in the mean and uniformly on compact subsets,
presents  the rescaled $($according to $(\ref{RR1}))$ solution of
the LSE $(\ref{1})$ for any data $u_0 \in L^2_{\rho^*}(\ren)$; and


\end{lemma}

\noindent{\em Proof.} (i) \underline{\em  Construction of
polynomial eigenfunctions}. Of course, the necessary discrete
spectrum \ef{SpecN} follows from \ef{PP1}. We now intend to show
how to obtain these results directly from the differential
operator.

 Thus, $\psi^*(y)$ is a polynomial. If its
degree is $k$, then
 $$
  \tex{
  \psi^*(y) = \sum\limits_{j=0}^s P_j(y),
  }
 $$
 where $P_j(y)$ is a homogeneous polynomial of degree $k-2m j$ with
 $s=\big[\frac k{2m}\big]$, denoting the integer part. Since by
 the Euler identity
 $$
  \tex{
 - \frac 1{2m} \sum\limits _{j=1}^N y_j \frac{\partial P_0(y)}{\partial y_j} = - \frac k{2m} P_0(y)
 = \l P_0(y),
  }
  $$
we see that $\l = -\frac k{2m}$ and $P_0(y)$ may be an arbitrary
homogeneous  polynomial of degree $k$. Other polynomials $P_j(y)$
are then defined as follows:
 $$
  \tex{
 P_j(y) = \frac 1{j !} \big(\ii(-\D)^{m}\big)^j P_0(y), \quad j=1,...,s.
  }
 $$
 This structure of $\psi^*(y)$ implies the completeness of the set of eigenfunctions in
 $L^2_{\rho}(\ren)$. In the second-order case $m=1$, this construction leads to the
 generalized Hermite polynomials, which were already introduced in \ef{kam2}.
Note that the polynomial structure of adjoint eigenfunctions
follows from the expansion (\ref{w1New}), where the coefficients
of initial data $\tilde M_\b(w_0)$ in (\ref{Mmom}) are the dual
products of $w_0 \in L^2_{\rho^*}(\ren)$ and $\tilde \psi^*_\b
\in L^2_{\rho}(\ren)$. This implies that each $\psi^*_\b(y)$
 is a finite linear combination of
elementary polynomials $y^\g$.

  We now fix $P_0(y) = y^\b$,  so that, for (extended)  eigenfunctions
 $\{\psi_\b\}$ of ${\mathbb B}$ in (\ref{eigen}), the corresponding adjoint eigenfunctions
 take the form
 \begin{equation}
 \label{psi**1}
 \fbox{$
  \tex{
 \psi_\b^*(y) = \frac 1{\sqrt{\beta !}}
 \Big[ y^\b + \sum\limits_{j=1}^{[\frac{|\b|}{2m}]} \frac 1{j !}\big(\ii(-\Delta)^{m}\big)^j y^\b \Big].
  }
  $}
 \end{equation}
 We also call \ef{psi**1} the {\em generalized Hermite polynomials}. For
 $m=1$, up to normalization constants, these coincide with those
 given by the classic generating formula \ef{kam2}.

\ssk

\noi (ii)  \underline{\em Completeness and closure}. This is the
well-known fact that polynomials $\{y^\beta\}$, which are
higher-order terms in any eigenfunction $\psi_\b^*$, are complete
in suitable weighted $L^p$-spaces; see \cite[p.~431]{KolF}.
Closure is associated with the eigenfunction expansion \ef{PP1}.






\noi (iii) \underline{\em  Evolution closure} follows from
\ef{PP1}.  $\qed$



\section{{\bf Application I:} Evolution completeness of $\Phi $ in $\LLL$, sharp estimates in $\re^{N+1}_+$,
 and some extensions}
 \label{S5Compl}

\subsection{Linear PDEs}
Our first result is about the following refined asymptotic
scattering:

\begin{theorem}
\label{Th.1} Consider the Cauchy problem $(\ref{1})$ for $u_0 \in
\LLL$ and $u_0 \not = 0$. Then there exists a finite $l \ge 0$
 and a function $\varphi_l(y)$, such that, as $t \to +\iy$,
 \beq
 \label{tt1}
  \tex{
  u(x,t)= t^{-\frac {N+l}{2m}}\big[ \varphi_l\big(\frac
  x{t^{1/2m}}\big) +O\big(t^{-\frac 1{2m}}\big)\big]
   }
   \eeq
   uniformly on compact sets in $y=\frac
  x{t^{1/2m}}$,
   where $\varphi_l(y)$ is a nontrivial superposition of extended eigenfunctions
 $\{\psi_\b, \, |\b|=l\}$  of ${\mathbb B}$ from the corresponding
 finite-dimensional eigenspace.

\end{theorem}

Of course, this is a corollary of our convergence analysis of the
series \ef{gsol} and \ef{w1New}, where $l$ is the minimal
multiindex  length $|\b|=l$, for which $\tilde M_\b(w_0) \not =
0$.
 As in Agmon's classic results for the parabolic case
(see zero set applications of advanced Agmon--Ogawa estimates in
\cite{XYCh} for parabolic PDEs for $m=1$), a super-fast decay in
\ef{tt1} corresponding to $l=\iy$ implies $u(x,t) \equiv 0$, so
that (other topologies are meant)
 \beq
 \label{tt2}
 |u(x,t)| \le t^{-K}\,\,\, \mbox{as $t \to \iy$ for any $K \gg 1$}
 \LongA u=0.
  \eeq

Further extensions of the above classification of $t \to \iy$
behaviour are rather straightforward for asymptotically small
perturbations of the LSE such as
 \beq
 \label{tt3}
  \tex{
 u_t= -\ii (-\D)^m u + \sum\limits_{(0 \le |\g|<2m)} a_\g(x,t) D^\g u
 \inB \ren \times \re_+,
  }
  \eeq
with, say, bounded complex-valued coefficients $\{a_\g\}$, which
 decay sufficiently fast:
 \beq
 \label{tt31}
 |a_\g(x,t)|=o\big(t^{\frac{|\g|}{2m}-1}\big) \asA t \to +\iy
 \,\,\, \mbox{for any}
 \,\,\, 0 \le |\g| \le 2m-1.
  \eeq
  We then deal with solutions $u(\cdot,t) \in \tilde L^2_\rho(\ren)$, for
  which eigenfunction expansions make sense.
 Then after scaling \ef{Resc11N}, we arrive at the same equation
  \ef{weq}
 with  asymptotically (and exponentially if necessary) small perturbations, which can be tackled
 by the eigenfunction expansion techniques, though some parts of
 the study can be indeed technical.
 Some extra efforts are necessary to tackle convergence properties of such
 series that look rather technical, though can be involved in some places.


\subsection{NLSE: on a ``centre subspace" behaviour}
 \label{S8.2}

There are no doubts that, at least partially,  the classification
results for decaying as $t \to \iy$ solutions can be extended to
the NLSE \ef{6}. Indeed, the nonlinear term then has the form
 \beq
 \label{tt4}
\ii |u|^{p-1}u \equiv a_0(x,t) u \whereA a_0(x,t)=\ii
|u(x,t)|^{p-1} \to 0 \,\,\,
 \mbox{as} \,\,\, t \to +\iy
 \eeq
 on such small solutions. Therefore,
 one can expect that
 the perturbation techniques are also effective here.

 Agreeing with that and do not performing this routine, but sometimes
  technical  analysis, we would
 like to stress our attention to the principal fact showing that
 even this looking purely  perturbation approach is not
 straightforward. Namely, we next formally show that the NLSE
 \ef{6} can admit some small solutions with a complicated asymptotics
 corresponding to {\em centre subspaces} of the rescaled operators.

  Bearing in mind the typical
 linear behaviour \ef{tt1}, where $l=0,1,2,...$\,, we perform
 in (\ref{6}) the standard
 rescaling as in \ef{Resc11N},
 \beq
 \label{tt5}
  \tex{
 u(x,t) = (1+t)^{-\frac{N+l}{2m}} v(y,\t), \quad y = \frac x{(1+t)^{1/{2m}}},
 \quad \t = \ln(1+t).
  }
 \eeq
 The rescaled solution $v(y,\t)$ satisfies the perturbed
 equation
  \begin{equation}
  \label{eqPer}
   \tex{
  v_\t = \big({\bf B}+ \mbox{$\frac l{2m}$} \, I\big) v  + \ii  {\mathrm e}^{-\g_l \t}
  |v|^{p-1}v \whereA
\g_l = \frac{(p-1)(N+l)}{2m} - 1.
   }
   \eeq

  It follows that there exists a sequence of critical exponents
  $\{p_l, \, l=0,1,2,...\}$ such that
  \beq
   \label{cr1}
    \tex{
  \g_l = 0 \quad \mbox{iff} \quad p=p_l=1+ \frac{2m}{N+l}, \quad l
  \ge 0.
   }
   \eeq
   In these critical cases, \ef{eqPer} yields  the autonomous equation
\begin{equation}
  \label{eqPer1}
  v_\t = ({\bf B}+ \mbox{$\frac l{2m}$}\, I) v  + \ii |v|^{p-1}v.
  \end{equation}
  Then,
we are looking for solution $v(\cdot,\t)$ with the behaviour for
$\t \gg 1$ close to the centre subspace of the linearized operator
${\bf B}+\frac l{2m}\,I$. Such a centre subspace asymptotic
dominance assumes that in the eigenfunction expansion
 of the solution
  \beq
  \label{ee3}
   \tex{
  v(y,\t)= \sum_{(\b)} c_\b(\t) \psi_\b(y)
 }
  \eeq
 the leading term for $\t \gg 1$ corresponds to an
eigenfunction $\varphi_l$ of ${\bf B}$,
 $$
  \tex{
 \big({\bf B}+ \mbox{$\frac l{2m}$} \,I \big) \varphi_l = 0,
 }
 $$
 i.e., $\varphi_l$ belongs to the centre subspace of the linearized operator
 ${\bf B}+ \frac l{2m}\,I$ in the perturbed equation (\ref{eqPer1}).
Hence,
 we suppose that
 \begin{equation}
 \label{vexp22}
  v(\t) = a_l(\t) [\varphi_l + o(1)] \quad \mbox{as} \,\,\, \t \to
 +\infty
 \end{equation}
 and this asymptotic equality can be differentiated in $y$ and
 $\t$.
 Then  the  equation for the leading expansion
 coefficient $a_l(\t)$ takes the form (recall the convention
 \ef{psiNN} for the metric)
 \begin{equation}
 \label{aleq}
 \dot a_l =  \ii |a_l|^{p_l-1}a_l  [c_l + o(1)] \quad \mbox{for} \,\,\, \t \gg 1, \quad
 \mbox{where} \,\,\, c_l=
 \langle |\varphi_l|^{p_l-1} \varphi_l, \varphi_l^* \rangle,
 \end{equation}
\noi and $c_l \not = 0$ as should be  assumed.
 Here, for the first (and the last) time, we need the values of  generalized  {\em
 extended
linear functional} to be introduced in Appendix D. We must admit
that
 it is not easy to evaluate such coefficients $c_l$, even numerically in 1D.
Moreover, we still do not know any efficient method to get these
values, so we will need further speculations.

However, to confirm that such a formal analysis actually makes
sense, we present a simple {\em explicit} example of such
solutions on a centre subspace (in fact, on a manifold).

\ssk

\noi{\bf Example: explicit  centre subspace periodic solutions for
$m=1$, $l=0$.} Obviously, the
 simpler case is $l=0$ and $m=1$,
where the eigenspace is 1D, so that $\varphi_0=F$ and hence, by
\ef{Ga1},
 \beq
 \label{kk1}
  \tex{
 |F(y)|=
 b_1= \frac 1{(4 \pi)^{N/2}}
 \LongA
 c_0= \langle
 |F|^{p_0-1} F, \psi_0^* \rangle= b_1^{\frac{2m}N},
 }
  \eeq
  so that $c_0$ is real.
  We have used the convention $\int F=1$, which, for $N=1$, holds
  in the usual sense, since the improper integral converges.
Then substituting into \ef{eqPer1}, with $l=0$, yields the
following explicit solution:
 \beq
 \label{sol12}
 \tex{
 v(y,\t)=a_0(\t) F(y) \LongA \dot a_0 = \ii \,|a_0|^{\frac 2N} a_0
 \LongA  \frac{\mathrm d}{{\mathrm d}\t} \, |a_0(\t)|^2=0,
 }
  \eeq
  i.e., this explicit centre subspace behaviour is a {\em periodic
  orbit}. Note that this is not an $L^2(\ren)$ or any $L^p(\ren)$
  solution, since $F(y) \equiv \psi)(y)$ is the first (extended)
  eigenfunction of the operator ${\mathbb B}$, and also a standard
  eigenfunction of $\BB$ in the extended space of eigenfunction
  expansion closures; see Appendix D. Therefore, we do not know
  whether a centre subspace solution \ef{sol12} can have a stable
  orbit connection with more customary $L^2$ and other solutions of Schr\"odinger equation theory.

  \ssk

\noi{\bf   For general $ m \ge 2$ and arbitrary $l \ge 0$},
 we need further arguments concerning the coefficients $c_l$ in
 \ef{aleq}. For instance,
   using the
  continuity argument with respect to $p$,
  we again observe ``almost" real values, since:
 \beq
 \label{hh1}
c_l = 1 \,\,\, \mbox{by orthonormality, if $p_l=1$, so  $c_l
\approx 1$ if $p_l \approx 1$ for all $l \gg 1$}.
 \eeq
 However, such an
asymptotic result assumes a technical proof, which falls out of
the scope of the present analysis. Overall, for real coefficients
$c_l$, the system (\ref{aleq}) is Hamiltonian,
 $$
 \tex{
 \frac{\mathrm d}{{\mathrm d}\t} \, |a_l(\t)|^2=0,
 }
 $$
 so, as above for $l=0$, it describes a {\em periodic behaviour}
 close to the centre subspace.

  For complex valued $c_l$, such a
 centre subspace behaviour can be more complicated.
 For instance, for $c_l=\ii \hat c_l$, with $\hat c_l>0$
  (a rather hypothetical situation to be used as an
  illustration only),
 integrating
\ef{aleq}
  yields the following
 decaying functions:
 \begin{equation}
 \label{albeh1}
 a_l(\t) =  C_l \t^{-\frac 1{p_l-1}} (1+o(1)) \quad \mbox{as} \,\,\,
 \t \to +\infty,
 \end{equation}
 where $C_l \in {\mathbb C}$ is a constant.
In terms of the original $(x,t,u)$-variables, such a behaviour
takes a form of a logarithmically perturbed linearized pattern
 \begin{equation}
 \label{ucr33}
  \tex{
 u(x,t) =C_l (t \ln t)^{-\frac{N+l}{2m}}  \big[\varphi_l\big(\frac x {t^{1/{2m}}}\big) + o(1)
 \big]
 \quad \mbox{as} \,\,\, t \to +\infty.
 }
 \end{equation}

As we have mentioned,  similar asymptotic patterns can be
constructed for the ``stable", {\em defocusing NLSE}
 \beq
   \label{6dd}
   u_t= -\ii (-\D)^m u - \ii |u|^{p-1} u \inB \ren \times \re
   \quad
   (p>1);
    \eeq
 see
 \cite{Ken06, Tao07, Vis07, Wang08} for key references and results
 concerning
 \ef{6dd} for $m=1$, as well as recent papers \cite{Bar10, Bar10Ring, Guo10, Miao09, Paus09, Paus10, Zhu10}
  (and references/short surveys
  therein) for $m=2$, i.e., for the {\em biharmonic} nonlinear Schr\"odinder  equation
   \beq
   \label{Sht4}
   \ii u_t+ \D^2 u= \pm |u|^{p-1} u \inB \ren \times \re.
    \eeq
For \ef{6dd}, such a centre subspace approach looks  even more promising
than for the unstable PDE \ef{6} admitting also blow-up in these
ranges.
 More flexibility is added when replacing the nonlinear term by a
 more general one,
  $$
  \pm \ii |u|^{p-1} u \mapsto d \,|u|^{p-1} u \whereA d = a + \ii
  b \in {\mathbb C}, \,\,\, ab \not =0.
  $$

 A rigorous
justification of the centre manifold-like patterns of a periodic
or \ef{ucr33}-type is a difficult open problem, which we do not
touch here. Notice that even existence of an invariant manifold
(in which functional setting?--in an extended space of closures as
in Appendix D?) is a very difficult problem for such Hermitian
spectral theory dealing with the pair $\{{\mathbb B},\,{\mathbb
B}^*\}$.
 This is regardless the good spectral properties of ${\mathbb B}$
listed in Lemma \ref{lemspec} and also its sectorial setting in
the topology of $l^2_\rho$ in Proposition \ref{PrDens}, which
however suggest a certain confidence that this behaviour can be
verified by using the powerful machinery  of classic invariant
manifold theory, \cite{Lun}.


\section{{\bf Applications II and III:} Local structure of nodal sets and unique continuation}
 \label{S6}

  \subsection{Application II:  blow-up formation of multiple zeros for linear PDEs (a Sturmian theory)}

Next, we arrive at the following classification of  zeros of
solutions of the LSE (\ref{1}):

\begin{theorem}
\label{Th.2}
 Consider the Cauchy problem $(\ref{1})$ for $u_0 \in
\LLL$ and $u_0 \not = 0$. Assume that the corresponding solution
$u(x,t)$
 creates a zero at a point $(0,T)$, i.e.,
$u(0,T)=0$.
 Then there exists a finite $l \ge 1$
 and a generalized Hermite polynomial $\varphi_\b^*(y)$ such that
 \beq
 \label{tt1n}
  \tex{
  u(x,t)= (T-t)^{\frac {l}{2m}} \big[\varphi_l^*\big(\frac
  x{(T-t)^{1/2m}}\big)+o(1)\big] \asA t \to T^-
   }
   \eeq
   uniformly on compact sets in $y=\frac
  x{(T-t)^{1/2m}}$,
   where $\varphi_l^*(y)\not \equiv 0$ is a superposition of
   polynomial (extended) eigenfunctions
 $\{\psi_\b^*, \, |\b|=l\}$  of ${\mathbb B}^*$ from the corresponding
 eigenspace ${\rm ker}\,\big({\mathbb B}^*-\frac l{2m}\big)$.

\end{theorem}

Since $\varphi_l(y)$ is  a generalized Hermite polynomial, the
multiple zero of $\RR u(x,t)$ (or, equivalently,  $\II u(x,t)$)
occurs    at the point $(0,T^-)$ by ``blow-up focusing" of several
zero-surfaces $\{x_\g(t)\}$ of $  \RR \varphi_l^*\big(\frac
  x{(T-t)^{1/2m}}\big)$, which move according to the scaling
  blow-up law
   \beq
   \label{ll1}
   x_\g(t) \sim (T-t)^{\frac 1{2m}} \to 0 \asA t \to T^- \quad (|\g|\le l).
 \eeq

The result \ef{tt1n} follows from the series \ef{PP1}, for which
 the rescaling \ef{RR1} is performed relative to the time moment $T$
 rather than 1.
 In a natural sense, the countable family of the types of asymptotics \ef{tt1n}
describes the sharp ``micro-turbulent" structure of the PDE
\ef{1}, since, by evolution completeness, on smaller space-time
scales, the solution behaviour is trivial (a constant one mostly).
 In other words, \ef{tt1n} exhausts all possible micro
 configurations that can be created by the LSE (and also by many
 other related semilinear and quasilinear PDEs admitting similar
 blow-up rescaling and Hermitian spectral properties; see below).
 See also
   \cite{2mSturm} for parabolic and other
 real-valued PDEs.

\subsection{Application III: Unique continuation}

Various
 classic and other well-known new  unique continuation results for linear and nonlinear
 Schr\"odinger-type PDEs can be found in \cite{Dong07, Ion04, Ion06, Ken02}, where
 further references are available. These directions on uniqueness PDE
 theory have their origins in many principal works in the twentieth century including such
 a classic path as Holmgren
 (the starting point, 1901)--Carleman (1933)--Myshkis (1948)--Plis (1954)-- Calderon (1958)--Agmon--Nirenberg...\,.
 Here we present an example of a slightly different
  type of a ``blow-up micro-scale
 uniqueness" study based on the spectral properties of ${\mathbb B}^*$,
 which is responsible for blow-up scaling of the PDE.

 Thus, all the types of nodal sets of
  zeros for the LSE \ef{1} are exhausted by the zero structures of
  $\RR$ or $\II$ of
  all the generalized Hermite polynomials $\Phi^* =\{\psi_\b^*(y)\}$ given
  by \ef{psi**1} including arbitrary linear combinations on all
  the eigenspaces. We fix this in the following rather
  unusual unique continuation theorem:

  \begin{corollary}
   \label{Cor.Un}
   Let under the hypotheses of Theorem $\ref{Th.2}$, the nodal set
   of the real part $($or, equally, of the imaginary part$)$ of the solution
    \beq
    \label{NN1}
    {\mathcal N}(u)=\{x \in \ren, \,\, t\in \re: \,\,\, \RR u(x,t)=0\}
 \eeq
 has a nontrivial component that evolves as $t \to T^-$ in a
 manner that does not asymptotically  match the zero sets of any finite linear combinations of the real parts of the
  generalized
 Hermite polynomials from $\Phi^* $. Then
 \beq
 \label{NN2}
  u(x,t) \equiv 0 \inB \ren \times \re.
   \eeq
   \end{corollary}

Of course, this is just a conventional version of the uniqueness
result that is based on the eigenfunction expansion \ef{PP1},
which can be expressed in a different and more reliable manner.
 For instance, e.g., depending on $l$ and also much on $m \ge 2$,
  for the case, where the resulting polynomial $\RR \varphi_l^*(y)$
 does not change sign at all (so that \ef{NN1} is locally empty), we then have to postulate
 just the asymptotic behaviour such as \ef{tt1n}. Then the alternative (and more correct
 and universal) sounding of the unique continuation result will be as follows:
 \beq
 \label{NN3}
 \mbox{if $u(x,t)$ violates any of non-trivial asymptotics (\ref{tt1n}) near zero, then $u \equiv 0$.}
  \eeq

   Recall that, typically, for the real-valued evolution linear or nonlinear PDEs
   with interior regularity, unique continuation theorems stated in
   the pointwise sense deals with zeros of infinite order in the
   following manner: if, in a natural integral mean sense,
    \beq
    \label{NN31}
   \mbox{$u$ has an infinite-order zero at $(0,T)$, then
$u(x,t) \equiv 0$};
   \eeq
   see \cite{XYCh} and \cite{2mSturm} for further references and
   results for parabolic PDEs (such results are also known for the LSEs and are proved
   by iterating
   Carleman's classic estimates). Of course, \ef{NN31} becomes trivial
   for analytic solutions (though extensions to smooth
   non-analytic ones along the lines discussed below makes deep
   sense), so that we present a new pointwise uniqueness version \ef{NN1}, which
   looks not that trivial.

It is natural to expect that the above classification of all the
possible  zeros remains for the perturbed LSEs such as \ef{tt3},
with arbitrary bounded coefficients $\{a_\g\}$. We then need to
assume that $u(x,t)$ is locally good enough close to the point
$(0,T^-)$, and at least, $u(\cdot,t) \in \tilde L^2_{\rho}(\ren)$
(see details in Appendix D), so  we can use the corresponding
eigenfunction expansions endowed with a strong enough  topology of
convergence on compact subsets,
 which are
 sufficient to detect and identify the zero structure of
solutions. At least, we need convergence a.e., which is guaranteed
by the $L^2_{\rho}$-metric.
 Nevertheless,  the
pointwise sense of such expansions will possibly demand some extra
hypotheses that are not discussed here.
 It is known that, even in the parabolic case,
such a Sturmian theory on zero sets leads to a number of technical
difficulties; see \cite{2mSturm}, where further references are
given and   applications to other classes of PDEs are discussed.
 Note also that such extensions can be applied to related partial differential
 inequalities (PDIs), e.g., for
  \beq
  \label{PI1}
  |u_t + \ii (-\D)^m u | \le C(|u|+|\n u|+...+ |D^{2m-1} u|),
   \eeq
 where $C>0$ is a constant. One can see that the right-hand side is always negligible after rescaling
 \ef{RR1}, so it does not affect the asymptotic  zero
 classification \ef{tt1n}.


\subsection{The NLSE: similar local zero set behaviour}

 For the NLSE \ef{6}, the local  zero
evolution remains unchanged since no centre subspace patterns are
available.   Namely, assuming again that a zero occurs at $(0,T)$,
we perform for \ef{6} the scaling
 \beq
 \label{ss111}
  \tex{
 u(x,t)=v(y,\t), \,\,\, y= \frac x{(T-t)^{ 1/{2m}}}, \,\,\,
 \t=-\ln(T-t) \,\, \Longrightarrow \,\,
  v_\t = \BB^* v \pm  {\mathrm e}^{-\t}\, \ii |v|^{p-1} v.
 }
 \eeq
 In other words, the nonlinear term near the  zero always
 creates an exponentially small as $\t \to +\iy$ perturbation of the dynamical
 system for the expansion coefficients. Hence, it is very unlikely that
 this can somehow essentially affect the local zero structure near
 $(0,T)$. Recall again that a rigorous analysis is rather involved
 even in simpler parabolic cases, \cite{Eg4, 2mSturm}. Of course,
 close to the  zeros, the NLSE \ef{6} falls into the
 scope of the PDI \ef{PI1}.



\begin{appendix}
\section*{Appendix A. ``Schr\"odinger equations" are most popular
in the twenty first century according to {\tt MathSciNet}}
 \label{S3M}
 \setcounter{section}{1}
\setcounter{equation}{0}

\begin{small}

We have claimed in Section \ref{S1} that Schr\"odinger PDEs are
most popular currently in mathematical PDE theory. To ``prove"
this statement, we present below Table 1 on  citations of the
different types of equations in the {\tt MathSciNet}. One can get
these results by typing in the box $\fbox{$\mbox{Title}$}$ the
corresponding type of equations, i.e., ``hyperbolic (elliptic,
parabolic, etc.) equation". Each resulting page in {\tt
MathSciNet} contains a variable number of papers (in average, it
is about 20, but can be more than that), so for convenience we
count the number of pages. Naturally, at this moment, we put a
blind eye on the obvious fact that there can be, and indeed many,
papers on various PDEs, whose title do not contain those two
words. Anyway, we believe that this is a reliable statistics and
the table reflects some true information.


\begin{table}[h]
\caption{Citations in {\tt MathSciNet}, 8th September, 2008}  
\begin{tabular}{@{}lll}
 PDE type in the title & All matches & In 2007--08
 \\\hline
\ssk
  ``hyperbolic equation" & 3321 & $\sim 5$ pages
   \\
   \ssk
 ``elliptic equation" & 7118 & $\sim 19$ pages
   \\
   \ssk
    ``parabolic equation" & 6972 & $\sim 18$ pages
   \\
   \ssk
    ``Schr\"odinger equation" & 5264 & $\sim {\bf 25!}$ pages
 \\
\end{tabular}
\end{table}

\noi Thus, the absolute winners  for the total number of papers
(via ``All matches"), which have in the title ``... equation" are
``elliptic" and ``parabolic" ones (recall that the {\tt
MathSciNet}  operated that time with about 3.5 million papers,
written by about 200\,000 mathematicians). However, during two
years 2007--08, more papers were published for ``Schr\"odinger"
ones: listing all of them takes about the record {\bf  25!} pages.
This adds extra bits to our motivation of developing a refined
spectral theory for the Schr\"odinger PDE \ef{1}.


\end{small}
\end{appendix}

\begin{appendix}
\section*{Appendix B. {\bf Application IV:} Towards boundary point regularity theory}

 \label{SReg}
 \setcounter{section}{2}
\setcounter{equation}{0}

\begin{small}


We now present another application of a refined $\{{\mathbb
B},{\mathbb B}^*\}$-spectral theory to the problem of regularity
of {\em boundary characteristic points} for Schr\"odinger
equations such as \ef{2} and \ef{1}. These regularity issues were
always in the core of general potential theory, which in its turn
represents one of the most classic directions of linear and
nonlinear PDE theory initiated already by Dirichlet himself in
1820s. We refer to Maz'ya's monographs with collaborators
\cite{KMR1, KMR2} for the history and key results on elliptic
PDEs, as well as to   recent surveys in \cite {GalPet2m,
GalMazf(u)} devoted more to parabolic PDEs, whose approaches and
results will be essentially used later on.

Of course, there are many very strong and classic boundary
regularity results for  Scr\"odinger equations \ef{2} and \ef{1},
which are explained in several key papers mentioned in the
introduction.
 In particular, Schr\"odinger (or similar and often equivalent beam-type) equations in {\em non-cylindrical domains}
 have been studied in the 1960s by J.-L. Lions and E.~Magenes, and by G.A.~Pozzi, and later on by T.~Gazenave and others. We refer to papers \cite{Ant10, Bern01, Clark07}, where further references and results can be found.
The existence-regularity results therein in principal cannot treat {\em characteristic} boundary points, where, as we show, the continuity of the solutions is rather tricky.
We hope that our brief regularity
exposition based on spectral theory of the pair $\{{\mathbb
B},{\mathbb B}^*\}$ will insert some  new features
 to this classic area, which were not observed before.

\subsection{Regular boundary points and general asymptotic problem}

Without loss of generality, in order to explain the main
ingredients of the boundary point regularity,  we consider the
simplest case $m=N=1$, i.e., the 1D  {\em second-order
Schr\"odinger equation}:
 \beq
 \label{R.1}
 u_t = \ii u_{xx} \quad \mbox{in} \quad Q_0,
  \eeq
where $Q_0$ is a typical domain, for which $(0,0)$ is a {\em
characteristic point}, i.e., the straight line  $\{t=0\}$ touches
the lateral boundary of $Q_0$. We define the {\em backward
parabola}\footnote{In $\ren$, it is a {\em backward paraboloid}
with a quite similar study as for $N=1$, though, of course, it
becomes more involved; see those typical boundary regularity
features in \cite{GalMazNSE} to 3D Navier--Stokes equations.}
$Q_0$ as follows, assuming that it has a single finite right-hand
lateral boundary (we will see why this is necessary):
 \beq
 \label{R.2}
 Q_0=\{-\iy < x< R(t),
 \,\,\, t \in (-1,0)\},
  \eeq
  where $R(t)>0$  for all $t \in [-1,0)$ is a given sufficiently smooth
  function satisfying
   \beq
   \label{R.3}
    R(t)  \to 0^+ \asA t \to 0^-.
     \eeq
    Finally, we pose the Dirichlet boundary
     condition:
     \beq
     \label{R.4}
      u(R(t),t)=0 \forA t \in (-1,0),
       \eeq
       and prescribe bounded smooth $L^2$-initial data at $t=-1$ (without loss of generality,
       we allow $u_0(x)$ too decay exponentially fast as
       $x \to -\iy$ or to even to be compactly supported):
 \beq
 \label{R.5}
 u(x,-1)=u_0(x) \forA x \in (-\iy, R(-1)).
  \eeq
Overall, \ef{R.1}--\ef{R.4} is a well-posed initial-boundary value
problem, and we assume that it possesses  a classic bounded
solution up to the characteristic moment $t=0^-$.

Thus, the point $(0,0)$ is called {\em regular} in Wiener's
classic sense if
 \beq
 \label{R.6}
 u(0,0)=0
  \eeq
  for any such data $u_0$,
  i.e., there is continuity along the boundary,
  and {\em irregular} otherwise (if \ef{R.6} fails for some $u_0$). Our main goal is to show how to
  answer the following question:
   \beq
   \label{R.7}
    \fbox{$
   \mbox{for which lateral boundaries given by $R(t)$, $(0,0)$ is regular (irregular).}
   $}
   \eeq

In fact, this follows the canonical  regularity statement by
I.G.~Petrovskii in 1934--35 \cite{Pet34Cr, Pet35}, who almost
completely solved the boundary regularity problem for the {\em
heat equation}
 \beq
 \label{R.8}
 u_t=u_{xx}.
  \eeq
  This led to his famous ``$\log\log$ backward parabola", meaning, in
  particular, the following remarkable and delicate results:
  \beq
 \label{PP1NN}
 \fbox{$
 \begin{matrix}
{\rm (i)}\,R(t)=2 \sqrt{-t} \,\,\sqrt{\ln|\ln (-t)|} \LongA
\mbox{$(0,0)$ is regular}, \andA \qquad\quad \ssk\\ {\rm (ii)}\,
R(t)=2(1+\e)\sqrt{-t}\,\, \sqrt{ \ln|\ln (-t)|}, \,\, \e>0 \LongA
 \mbox{$(0,0)$ is irregular}.
  \end{matrix}
  $}
  \eeq
In what follows, unlike many strong well-known results and
approaches, we follow an asymptotic {\em matching blow-up}
approach
to the present regularity problem, which was  developed in
\cite{GalPet2m} and \cite{GalMazf(u)} for  higher-order linear and
nonlinear parabolic PDEs, respectively.

Actually, of course, we are solving a more general problem on the
asymptotic behaviour of solutions $u(x,t)$ as $x \to 0$ and $t \to
0^-$, so our goal is as follows:
 \beq
 \label{PP12}
  \fbox{$
 \mbox{to describe sharp asymptotics of $u(x,t)$ at the blow-up
 point $(0,0^-)$,}
 $}
 \eeq
and essentially detect their dependence on the function $R(t)$
defining the right-hand lateral boundary. As usual, once \ef{PP12}
has been solved, one can check the regularity property \ef{R.7}.

\subsection{A comment: two boundary conditions must be non-Hamiltonian}

Let us comment on the non-symmetric shape of the domain in
\ef{R.2}. Using the symmetric shape (as can be done
 for the heat equation \ef{R.8} or for the
bi-harmonic one \ef{RO4})
 \beq
 \label{R.9}
\hat  Q_0=\{-R(t)<x< R(t),
 \,\,\, t \in (-1,0)\}
  \eeq
is not possible. Indeed, since \ef{R.1} with the Dirichlet
conditions on the lateral boundary of $\hat Q_0$  is a Hamiltonian
system with the $L^2$-conservation, as $t \to 0^-$, one observes
the concentration of this $L^2$-energy onto shrinking-to-zero
$x$-intervals, so that, obviously,
 \beq
 \label{R.10}
 \mbox{for domains (\ref{R.9}) and conservative boundary conditions,  $(0,0)$ is irregular,}
 \eeq
 that eliminated the regularity issue at all.
On the other hand, the symmetric domains such as \ef{R.9} are
admitted if the boundary conditions on the lateral boundary
violate the Hamiltonian (symplectic) $L^2$-conservation property.
It is not that easy to find such conditions for the second-order
equation \ef{R.1}. For instance, these could be the {\em Robin
ones} (a third kind condition) at $x=-R(t)$,
 $$
 u_x + \s u=0 \quad \mbox{for some $\s \in {\mathbb C}$},
  $$
  and then the regularity would mean the
continuity at $(0,0)$ along the right-hand boundary.

For the fourth-order LSE,
 \beq
 \label{R.11}
 u_t=- \ii u_{xxxx},
 \eeq
 a domain  \ef{R.9} also requires non-Hamiltonian lateral
 boundary conditions. By the identity
  \beq
  \label{Id1}
 \tex{
 \frac {\mathrm d}{{\mathrm d}t}\, \int |u|^2=\ii\,\Big(
 \big[\bar u_{xxx} u - u_{xxx} \bar u\big]^{R(t)}_{-R(t)}+
 \big[u_{xx} \bar u_x - \bar u_{xx} u_x \big]^{R(t)}_{-R(t)} \Big),
  }
  \eeq
one can see that the homogeneous Dirichlet conditions
 $$
 u=u_x=0 \quad \mbox{at} \,\,\, x=\pm R(t), \,\,\, t \in
  (-1,0),
 $$
 are indeed Hamiltonian, the $L^2$-norm of $u(\cdot,t)$ is
 preserved in such  domains shrinking to a point, so $(0,0)$ is always irregular.
 The same is true for
 the Navier-type conditions
  $$
  u=u_{xx}=0\quad \mbox{at} \,\,\, x=\pm R(t).
  $$
  Again, by \ef{Id1}, there exists the $L^2$-conservation, so the vertex regularity
  problem makes no sense and $(0,0)$ is irregular for any nontrivial initial data.

  On the other hand,
 the following conditions
  \beq
  \label{R.12}
  u=u_{xxx}=0
  \quad \mbox{at} \,\,\, x=\pm R(t),
   \eeq
in general, violate the
  $L^2$-conservation on shrinking domains a $t \to 0^-$, so the
  characteristic boundary regularity problem makes sense. We can
  apply our ``blow-up" scaling-matching approach to study
  regularity of the vertex for conditions like in \ef{R.12} or
  others of higher-order Robin-kind (but this will require a
  different boundary layer theory, see below).

The boundary point regularity approach  proposed here also covers
the problem \ef{R.11}, \ef{R.12} and other $2m$th-order ones with
various (non-symplectic) boundary data.

\subsection{Introducing slow growing factor $\var(\t)$}

Thus, we return to the canonical (and indeed looking very simple)
LSE--2 \ef{R.1} in the one-sided domain \ef{R.2}, with all
conditions already specified.
Then, similar to \ef{PP1NN}, we  introduce a one-sided {backward
parabola} at $(0,0)$ given by the function
 \beq
 \label{ph1}
 R(t)= \sqrt{-t}\,\, \var(\t) \whereA \t= - \ln(-t) \to + \iy
 \asA t \to 0^-.
  \eeq
In Petrovskii's criterion \ef{PP1NN},
 $$
 \var(\t) \sim 2 \sqrt {\ln \t} \asA \t \to +\iy,
 $$
 so that $\var(\t)$ is
 an unknown  slow growing function
 satisfying
 \beq
 \label{vv1}
  \tex{
 \var(\t) \to +\iy, \quad \var'(\t) \to 0, \andA \frac {\var'(\t)}{\var(\t)} \to 0 \asA \t \to
 +\iy.
  }
  \eeq
 Moreover, as a sharper characterization of the above class of
 {\em slow growing functions}, we use the following criterion:
  \beq
  \label{vv2}
   \tex{
   \big( \frac{\var(\t)}{\var'(\t)}\big)' \to \iy \asA \t \to +\iy
    \quad (\var'(\t) \not = 0).
    }
    \eeq
    This is a typical condition in blow-up analysis distinguishing classes of
    exponential (the limit is 0), power-like (a constant $\not =
    0$), and slow-growing functions. See \cite[pp.~390-400]{SGKM},
    where in Lemma 1 on p.~400, extra properties of slow-growing
    functions \ef{vv2} are proved to be used later on.

\subsection{First kernel scaling}

 By \ef{ph1}, we
perform the similarity scaling
 \beq
 \label{ph3}
  \tex{
  u(x,t)=  v(y,\t) \whereA y=\frac x{\sqrt{-t}}.
  }
   \eeq
   Then the rescaled function $v(y,\t)$ solves the
   rescaled equation
 \beq
 \label{ph4}
  \left\{
  \begin{matrix}
  v_\t= \BB^* v \equiv  \ii v_{yy}- \frac 12 \, y v_y
  \inB Q_0=\{-\iy <y < \var(\t), \,\,\, \t>0\}, \ssk\ssk\\
  v=0 \atA y= \var(\t), \,\,\t \ge 0,
  \qquad\qquad\qquad\qquad\qquad\qquad\qquad\quad
  \ssk\ssk\\
  v(0,y)=v_0(y) \equiv u_0(y) \onA
  (-\iy,R(-1)),\qquad\quad\qquad\qquad\quad\,\,\,
   \end{matrix}
   \right.
  \eeq
 where, by obvious (blow-up micro-scale) reasons, the rescaled
 differential expression
 $\BB^*$, defined as in \ef{RO1} for  $m=N=1$,
 appears.
 In view of the assumed divergence \ef{vv1}, it follows that our final analysis will
 essentially depend on the spectral properties of the corresponding restricted linear
 operator ${\mathbb B}^*$ on the whole  line $y \in \re$, i.e., we arrive at the necessity of
   Hermitian spectral theory developed above.
Note that, in particular, this  differs our regularity analysis
from several well-known ones such as  Kondrat'ev's classic results
of the 1966--67 \cite{Kond66, Kond67} (see a later survey
\cite{KondOl83}), where, as a rule, the rescaled boundary remains
asymptotically fixed, which  is key  for using the spectral
properties of the bundles of linear operators in locally compact
domains. We will present further more detailed comments on that
below.

 \subsection{Regularity of a fixed backward parabolae is not
 obvious}

First of all, we need to comment on the  regularity of the vertex
$(0,0)$ of the  {\em backward fundamental parabolae}:
 \beq
 \label{P1}
 R(t)=l\sqrt{-t}, \quad \mbox{i.e.,} \quad \var(\t) \equiv
 l={\rm const.}>0.
  \eeq
Then the problem \ef{ph4} is considered on the fixed unbounded
interval
 \beq
 \label{Ill1}
 I_l=\{-\iy<y<l\},
   \eeq
   so that the final conclusion entirely depends on
spectral properties of ${\mathbf B}^*$ in $I_l$ with Dirichlet
boundary conditions.
 Since we need a sharp bound on the first eigenvalue, the clear
conclusion on regilarity/irregularity becomes rather involved,
where numerics are necessary to fix final details. In addition,
 as we pointed out,  in more general setting for the
fundamental backward paraboloids  in $\ren$, the existence,
uniqueness, and regularity of solutions in Sobolev spaces  was
proved in a number of papers such as \cite{Mih61, Mih63I, Mih63II,
Fei71}, etc.
 Note that in \cite[p.~45]{Mih63I}, the zero boundary data were
 understood in the {\em mean sense} (i.e., in the $L^2$-sense along a sequence
 of  smooth internal contours
 ``converging" to the boundary).

Note that it is not obvious at all that the spectrum of $\BB^*$ on
$I_l$ (with a standard $L^2_\rho$-setting as $y \to -\iy$) is
real. However, one can expect that, by a continuity argument, the
``first" eigenvalue $\l_0=\l_0(l)$ depending on $l$ satisfies
 \beq
 \label{l11}
  \l_0(l)  \to 0 \asA \l \to +\iy.
  \eeq
  This reminds a standard asymptotic  fact from
classic perturbation theory of linear operators (see Kato
\cite{Kato}) that the spectrum of $\BB^*$ in $I_l$ approaches as
$l \to +\iy$ that in $L^2_\rho(\re)$, according to Lemma
\ref{lemSpec2}. Then along with \ef{l11}, one can also expect
convergence of the first eigenfunction:
 \beq
 \label{l12}
  \psi^*(y;l) \to \psi_0^*(y) \equiv 1 \asA l \to + \iy.
   \eeq
Notice that \ef{l12} contains some features of a {\em boundary
layer} that occurs as $l \to +\iy$, which we will be use in the
non-stationary limit $\var(\t) \to +\iy$ as $\t \to +\iy$.

Overall, the limit \ef{l11} reflects the possibility for $(0,0)$
to be regular or irregular for different values of $l>0$ depending
of the sign of ${\rm Re} \, \l_0(l)$ (or $\l_0(l)$ itself,
provided that it is real). In other words, we conclude as follows:
if the limit \ef{l11} is oscillatory, then the backward parabolae
\ef{P1} can be regular or irregular. Note that this happens for
the bi-harmonic equation \ef{RO4}, where the corresponding fixed
parabolae with
 $$
 R(t)=(-t)^{\frac 14}l
 $$
 is regular for $l=4$ but is irregular for $l=5$. These
 conclusions for \ef{RO4} were fully justified numerically only,
 \cite[\S~6]{GalPet2m}.

 Thus, the regularity analysis of the backward parabolae
 \ef{P1} for the operator \ef{R.1} remains open, and its complete
 solution cannot be done without using enhanced numerical methods.
Nevertheless, despite such a theoretical gap for constant $l$'s,
we proceed to study the regularity for unbounded functions
\ef{vv1}, which promises even greater mathematical challenge.

\subsection{The case $\varphi(\t) \to 0$ as $\t \to +\iy$ is
always regular}

Indeed, in this case, as follows from \ef{ph4}, $(0,0)$ is always
regular, since $\var(\t) \to 0$ and hence, by just the continuity
of the solution $v(y, \t)$, we have  $v(0,\t) \to 0$ as $\t \to + \iy$.
It is worth mentioning that this is a completely {\em rigorous} result.
Therefore, according to our blow-up asymptotic approach, the
regularity problem of the vertex $(0,0)$ is trivial. Note that it
is not that trivial via  general PDE approaches (not including blow-up
scalings); cf. much weaker assumptions in, e.g.,  \cite{Clark07}.

 \subsection{Second scaling: Boundary Layer (BL) structure}
   \label{S3.3}

   Meantime, we return to the case of unbounded functions $\var(\t)$'s.
   Then, using standard boundary layer concepts of Prandtl--Blasius
   developed in 1904--08, we observe that,
 sufficiently close to the right-hand lateral boundary of $Q_0$, it is
 natural to introduces the variables
  \beq
  \label{z1}
   \tex{
   z= \frac y{\var(\t)} \andA v(y,\t)=w(z,\t)
    \LongA w_\t=  \frac 1{\var^2} \, \ii w_{zz} - \frac 12\, z w_z +
    \frac {\var'}\var \, z w_z.
    }
    \eeq
We next introduce the standard BL-variables (the same as for the
heat equation)
  \beq
  \label{z2}
  \tex{
  \xi= \var^{2}(\t)(1-z) \equiv \var(y-\var), \quad \var^{2}(\t)\, {\mathrm d} \t={\mathrm d}s,
   \andA w(z,\t)= \rho(s) g(\xi,s),
  }
  \eeq
  where $\rho(s)$ is an unknown slow decaying (in the same natural sense,
 associated with \ef{vv2})
   time-factor depending on the
  function $\var(\t)$.
 We will use later on the fact that
  \beq
  \label{rhoa1}
  \rho=a_0(\t)(1+o(1)) \asA \t \to +\iy,
   \eeq
   where $a_0(\t)$ is the first Fourier coefficient of the
   solution $v(y,\t)$ relative to the adjoint basis $\{\psi_k^*\}$
   of the operator ${\mathbb B}^*$.

On substitution into the PDE in \ef{z1}, we obtain the following
perturbed equation (see details of a similar derivation in
\cite[\S~7.2]{GalPet2m}):
 \beq
 \label{z3}
  \tex{
  g_s= \AAA g - \frac 1{\var^2}\big( \frac 12+ \frac {\var'}{\var}\big) \xi g_\xi - \frac
  {\var'_\t}{\var}\, g_\xi,
\quad  \mbox{where} \quad  \AAA g= \ii g '' + \frac 12\, g'.
 } 
   \eeq
As usual in boundary layer theory,
we are looking for a generic
 pattern of the behaviour described by \ef{z3} on compact subsets
 near the lateral boundary,
  \beq
  \label{z4}
  |\xi| = o\big(\var^{-2}(\t)\big)
   \LongA |z-1| = o\big(\var^{-4}(\t)\big) \asA \t \to +
   \iy.
   \eeq
On these space-time compact subsets, the second term on the
right-hand side of  \ef{z3} becomes asymptotically small, while
all the others are much smaller in view of the slow growth/decay
assumptions such as \ef{vv2} for $\var(\t)$ and $\rho(s)$.

Then posing the asymptotic behaviour at infinity: $g(\xi,\t)$ is
bounded as $\xi \to +\iy$, and
 \beq
 \label{z5}
  \tex{
  \sup_{\xi}|g(\xi,s)| = 1
   \quad (\mbox{Hypothesis I for
  generic patterns}).
  }
  \eeq
This is a typical ``normalization by 1" condition from boundary
layer theory. Note that, in view of highly oscillatory nature of
any solutions of Schr\"odinger equations, we cannot normalize by 1
any component of $g(\xi,s)=h(\xi,s)+\ii w(\xi,s)$. The condition
\ef{z5} will be used for matching with the solution asymptotics in
the Inner Region.

  Then, as $\xi \to + \iy$,
   all the derivatives are assumed to be bounded, and we arrive at
  a standard stabilization issue of passing to the limit as $s \to + \iy$ in
   \ef{z3}, \ef{z5}. Assuming that, by the definition in \ef{z2},
  the rescaled orbit $\{g(s),\,\, s>0\}$ is uniformly bounded, by
  classic parabolic theory \cite{EidSys}, one can pass to the
  limit in \ef{z3} along a subsequence $\{s_k\} \to +\iy$. Namely,
  by the above, we have that, uniformly on compact subsets defined in
  \ef{z4}, as $k \to \iy$,
   \beq
   \label{z6}
   g(s_k+s) \to h(s) \whereA h_s=\AAA h, \quad h=0
   \,\,\,\mbox{at} \,\,\,\xi=0,
    \eeq
    and (cf. \ef{z5}) $h$ {\em is bounded} (i.e., being oscillatory) as $\xi
    \to + \iy$.
The {\em limit} (at $s=+\iy$) {\em equation} obtained from
\ef{z3}:
 \beq
  \label{z7}
   \tex{
  h_s= \AAA h \equiv \ii h_{\xi\xi}+ \frac 12\, h_\xi
  }
   \eeq
is a standard linear  PDE in the unbounded domain $\re_+$, though
it is governed by a non self-adjoint operator
 $\AAA$. We then need the following its property:  in an appropriate weighted
 $L^2$-space if necessary and under the hypothesis \ef{z5},
 the stabilization holds, i.e.,
  the $\o$-limit set of $\{h(s)\}$ consists of
 equilibria: as $s \to +\iy$,
  \beq
  \label{z10}
   \left\{
   \begin{matrix}
  h(\xi,s) \to g_0(\xi) \whereA \AAA g_0=0 \,\,\, \mbox{for}
  \,\,\,\xi>0,\qquad\qquad\qquad\quad \ssk\ssk \\
   g_0=0 \quad \mbox{for} \quad \xi=0, \quad
 \sup_\xi |g_0(\xi)| = 1 \,\,\,\mbox{(bounded in $\re_+$)}.\qquad\,\,
  \end{matrix}
  \right.
   \eeq
 The stationary problem in \ef{z10} can be easily solved to give
 the BL profile
 \beq
 \label{z12}
  \tex{
  g_0(\xi)= \frac 12\, \big(1-{\mathrm e}^{\ii \frac \xi 2}\big)
  \whereA
  |g_0(\xi)|=|\cos(\frac \xi 4)| \le 1.
  }
  \eeq

  We must admit that \ef{z10} and \ef{z12} {\em actually define} the
  class of solutions we are going to treat later on. Hopefully,
  this should be  a generic class. So,
we will not concentrate on the stabilization problem \ef{z10},
which reduces to a standard spectral study of $\AAA$ in a weighted
space. Actually, the convergence \ef{z6} and \ef{z10} for the
perturbed dynamical system \ef{z3} is the main {\sc Hypothesis
(H)}, which characterizes the class of generic patterns under
consideration, and then \ef{z5} is its partial consequence. Note
that the uniform stability of the stationary point $g_0$ in the
limit autonomous system \ef{z7} in a suitable metric will
guarantee that the asymptotically small perturbations do not
affect the omega-limit set; see \cite[Ch.~1]{AMGV}.
 Such a definition of generic patterns looks rather
non-constructive, which, however, is unavoidable  for higher-order
PDEs without positivity and order-preserving  issues.

\subsection{Inner Region expansion: towards regularity}
 \label{S3.4}

 We next proceed as in \cite[\S~7.3,\,7.7]{GalPet2m}.
 Namely,
in Inner Region, we deal with the original rescaled problem
 \ef{ph4}.

  In order to apply the standard eigenfunction expansion
 techniques by using the orthonormal set of polynomial
 eigenfunctions of ${\mathbb B}^*$ given in \ef{psi**1},
 as customary in classic PDE and potential theory, we extend
 $v(y,\t)$ by 0 for $y > \var(\t)$ by setting:
  \beq
  \label{a1}
 \hat v(y,\t) = v(y,\t)H(\var(\t)-y)=
  \left\{
   \begin{matrix}
   v(y,\t) \forA 0 \le y < \var(\t), \\
\,\, 0  \,\,\forA \,\,  y \ge  \var(\t),
\end{matrix}
 \right.
 \eeq
 where $H$ is the Heaviside function.
Since $v=0$ on the lateral boundary $\{y= \var(\t)\}$, in the
sense of distributions,
 \beq
 \label{a2}
 \hat v_\t= v_\t H, \quad \hat v_y= v_y H, \quad
 \hat v_{yy}= v_{yy} H-v_y\big|_{y=\var} \d(y-\var).
  \eeq
Therefore, $\hat v$ satisfies the following equation:
 \beq
 \label{a3}
  \hat v_\t= \BB^* \hat v + v_{y}\big|_{y=\var}\d(y-\var)
  \inB \re_+\times \re_+.
  \eeq
Since the extended solution orbit \ef{a1} is assumed to be
uniformly bounded in $L^2_{\rho^*}(\re)$, we use the converging in
the mean (and uniformly on compact subsets in $y$) the
eigenfunction expansion via the generalized Hermite polynomials
\ef{psi**1}:
 \beq
 \label{a4}
  \tex{
  \hat v(y,\t)= \sum_{(k \ge 0)} a_k(\t) \psi_k^*(y).
   }
   \eeq
 Obviously, this assumes the inclusion $v(\cdot,\t) \in \tilde
 L_\rho^2$, for all $\t \gg 1$,
 which, by classic regularity theory for Schr\"odinger equations,
 is not a restrictive assumption at all. Then,
   substituting \ef{a4} into \ef{a3} and using the orthonormality
   property
    \ef{bi1} in a ``{\em v.p.}" sense of (\ref{MEAN1})
      yields the following dynamical system for the
   expansion coefficients:
   \beq
   \label{a5}
   \tex{
   a_k'= \l_k a_k +  v_{y}\big|_{y=\var(\t)} \langle
   \d(y-\var(\t)), \psi_k \rangle
   \quad \mbox{for all $k=0,1,2,...\, ,$}
    }
    \eeq
    where $\l_k= -\frac k 2$ are real eigenvalues \ef{spec1}
    of ${\mathbb B}^*$. Recall that
    $\l_k<0$ for all $k \ge 1$. More importantly, the
    corresponding eigenfunctions $\psi_k(y)$ are unbounded and not
    monotone for $k \ge 1$ according to \ef{psi**1}. Therefore,
    regardless proper asymptotics given by \ef{a5}, these inner
    patterns cannot be matched with the BL-behaviour such as
    \ef{z5}, and demand other matching theory. Since these are not
    generic, the latter will be dropped, though can be taken into account
    for a full classification of (non-generic) asymptotics.


     Thus,
     one needs to
     concentrate on the ``maximal" first Fourier generic pattern associated with
 \beq
 \label{a6}
 k=0: \quad \l_0=0 \andA \psi_0^*(y) \equiv 1 \quad
 \big(\psi_0(y)=F(y)\big),
  \eeq
  where $F(y)$ is the ``Gaussian" profile \ef{Ga1}.
 Actually, this corresponds to a naturally understood ``centre subspace behaviour" for
 the equation \ef{a5}:
  \beq
  \label{a7}
  \hat v(y,\t) = a_0(\t)\cdot 1 + w^\bot(y,\t) \whereA w^\bot \in {\rm
  Span}\{\psi_k^*, \,\, k \ge 1\},
   \eeq
   and $w^\bot(y,\t)$ is then negligible relative to $a_0(\t)$.
   This is another characterization of our class of generic patterns,
   Hypothesis II.
   The equation for $a_0(\t)$ then takes the form:
\beq
   \label{a8}
   \tex{
   a_0'=   v_{y}\big|_{y=\var(\t)}  \psi_0(\var(\t)).
    }
    \eeq

We now   return to  BL theory established the boundary behaviour
\ef{z2} for $\t \gg 1$, which for convenience we state again:
 in the
rescaled sense, on the given compact subsets,
 \beq
 \label{9}
  \tex{
  v(y,\t) = \rho(s) g_0\big( \var^{2}(\t)(1- \frac
  y{\var(\t)}) \big)+...\,.
  }
  \eeq
  Of course, since the limit BL-profile \ef{z12} is uniformly bounded but essentially
  oscillatory, we are talking about matching of this
  BL-asymptotics with the constant one in \ef{a7} {\em only} in a
  natural
  ``average sense". E.g., after a standard integral averaging the oscillatory BL-profile, which eliminates
 the non-essential multiplier $\frac 12$ in \ef{z12},
 since\footnote{This is about averaging of $|g_0|^2$; averaging of
 $|g_0|$ will give  another (non-important) constant.}
  $$
  \tex{
  \frac 1L \, \int_{0}^L \cos^2(\frac \xi 4)\, {\mathrm d} \xi \to
  \frac 12 \asA L \to +\iy.
   }
   $$
Note that there are (finitely) oscillatory generalized Hermite
polynomials $\psi_k^*(y) \equiv H_k(\ii y)$ for large $k$ (see
\ef{psi**1}, $k=\b$ in 1D), but these all are {\em unbounded} as
$y \to \iy$, so cannot be matched with uniformly bounded
BL-expansions.

 By the matching of both Regions, one concludes that, for such
 generic patterns,
  \beq
  \label{10}
   \tex{
   \frac {a_0(\t)}{\rho(s)} \to 1 \asA \t \to +\iy
   \LongA \rho(s) = a_0(\t)(1+o(1)).
   }
   \eeq
 Then the convergence \ef{z6}, which by a standard regularity is assumed to be
  also true for the
 spatial derivatives, yields, in the natural rescaled sense, that,
 as $\t \to +\iy$,
  \beq
  \label{11}
   \tex{
 v_{y}\big|_{y=\var(\t)} \to \rho(s) \var(\t) \g_1 \to a_0(\t)
 \var(\t) \g_1 \whereA \g_1= -g_0'(0)= \frac \ii 2.
 }
   \eeq
 We again recall that such an estimate is assumed to be true for a
 fixed above generic class of solutions satisfying a proper
 stabilization property in the BL.

Thus, this leads to an asymptotic ODE for the first expansion
coefficient for generic
 patterns:
 \beq
 \label{12}
  \tex{
    {a_0'}= \hat \g_1 a_0 \var(\t) {\mathrm e}^{\ii \frac
    {|\var(\t)|^2}4}+... \whereA \hat \g_1= \frac \ii{2\sqrt{4\pi
    \ii}}= \frac {1-\ii}{4 \sqrt {2 \pi}}.
    }
    \eeq
This not-that-easy asymptotic ODE gives insight into main
difficulties that one can face while posing and studying the
problem on the boundary regularity of the vertex $(0,0)$ for
Schr\"odinger-type operators. To this end, we first derive the
real form of this system for
 \beq
 \label{ab1}
 a_0(\t)=b_0(\t)+ \ii d_0(\t),
 \eeq
 where these parts now satisfy the system:
  \beq
  \label{ab2}
  \left\{
   \begin{matrix}
   b_0\,'= \frac{\var(\t)}{4 \sqrt {2 \pi}} \big[(b_0+ d_0) \cos \big(\frac {\var^2(\t)}4 \big)
    + (b_0-d_0) \sin \big(\frac {\var^2(\t)}4 \big)\big]+...\,, \quad\ssk\\
d_0\,'= \frac{\var(\t)}{4 \sqrt {2 \pi}} \big[(-b_0+ d_0) \cos
\big(\frac {\var^2(\t)}4 \big)
    + (b_0+d_0) \sin \big(\frac {\var^2(\t)}4 \big)\big]+... \, .
 \end{matrix}
 \right.
  \eeq
Note that the regularity of $(0,0)$ assumes that both limits are
zero:
 \beq
 \label{ab3}
  b_0(\t) \andA d_0(\t) \to 0 \asA \t \to + \iy.
   \eeq
 The system \ef{ab2} shows a general range of problems that appear
  determining the conditions on the lateral boundary given by the
  function $\var(\t)$ to ensure \ef{ab3}. Of course, there is no
  any hope to guarantee \ef{ab3} via a kind of Osgood (1898)--Dini-like
  integral condition of Petrovskii's type for the heat equation
  \ef{R.8}. The latter one has a very simple form (see a derivation in the
  same lines of spectral properties of $\BB^*$ in \cite[\S~7.7]{GalPet2m}):
 \beq
   \label{993}
    \fbox{$
\mbox{$(0,0)$ is regular iff} \quad
 \int\limits^\iy \var(s)\, {\mathrm e}^{-
 \frac{\var^2(s)}4}\, {\mathrm d}s = +\iy,
  $}
  \eeq
 which was
   already
  obtained by Petrovskii in 1934 \cite{Pet34Cr, Pet35} (earlier related results
  were due to Khinchin, 1924, in a probability representation; see
  \cite[\S~3.2]{GalPet2m} for further details).

\ssk

As a clue to some hard features of the system \ef{ab2}, consider
as a ``toy-model" a single equation of a  similar form:
 \beq
 \label{ab5}
  \tex{
 b_0\,' \sim \frac{\var(\t)}{4 \sqrt {2 \pi}} b_0\cos \big(\frac {\var^2(\t)}4
 \big).
  }
  \eeq
Integrating one obtains that
 \beq
 \label{ab6}
  \tex{
  \ln |b_0(\t)| \sim
  \frac{1}{4 \sqrt {2 \pi}} \int\limits^\t \var(s)\cos \big(\frac {\var^2(s)}4
 \big)\, {\mathrm d}s.
  }
  \eeq
In particular, the regularity of $(0,0)$ demands that
 \beq
 \label{ab7}
  \tex{
\ln |b_0(\t)| \to - \iy \LongA \int\limits^\t \var(s)\cos
\big(\frac {\var^2(s)}4
 \big)\, {\mathrm d}s \to -\iy \asA \t \to + \iy.
 }
 \eeq
This immediately implies (this is true for the toy-model, but it
seems such a condition exists for the whole system \ef{ab2})
 \beq
 \label{ab8}
 \tex{
\int\limits^\t \var(s)\cos \big(\frac {\var^2(s)}4
 \big)\, {\mathrm d}s \,\,\, \mbox{converges} \LongA
 (0,0) \,\,\, \mbox{is irregular}.
  }
  \eeq

It is easy to see that for the power-type functions
 \beq
 \label{ab9}
  \var(\t) \sim \t^\a \LongA \a>1 \,\,\, \mbox{means integral
  convergence and hence
  irregularity}.
   \eeq
   Though $\var(\t)$ in \ef{ab9} is not a slow growing function,
   the above results are extended to that case. Hence, for any
   slow growing functions, the integral in \ef{ab8} always
   diverges.

\ssk

However, on the other hand, the divergence of the integral in
\ef{ab8}, e.g., for $\a \le 1$ in \ef{ab9} {\em does not imply any
regularity}. Indeed, according to \ef{ab7}, the divergence must be
to $-\iy$, so that to create a regular boundary point $(0,0)$ a
special {\em oscillatory cut-off of the boundary} is necessary.
This is such a ``correction" of the shape of the boundary that
eliminates the positive divergence part in the integrals in
\ef{ab7}. Such a procedure, which establishes a coherent behaviour
of the lateral boundary close to the characteristic point $(0,0)$
and oscillations of the rescaled kernel $F(y)$, is inevitable for
infinitely oscillatory fundamental solutions. This oscillatory
cut-off already appears for the bi-harmonic equation \ef{RO4}; see
\cite[\S~7.5]{GalPet2m}, where further related details can be
found.

\ssk

We hope that the above analysis correctly describes a full range
of problems that appear in the regularity analysis and also
correctly shows how  $\{{\mathbb B},{\mathbb B}^*\}$-spectral
theory naturally enters this classic PDE area.

\end{small}
\end{appendix}

\begin{appendix}
\section*{Appendix C. {\bf Application V:} Towards countable families of
 nonlinear eigenfunctions of the QLSE}
 \label{SBQ}
 \setcounter{section}{3}
\setcounter{equation}{0}

\begin{small}


\subsection{The QLSE and its applications}

This is a more striking and even controversial application of the
spectral $\{{\mathbb B},{\mathbb B}^*\}$ theory developed above.
Namely, we now turn to a $2m$th-order {\em quasilinear
Schr\"odinger equation} (the QLSE--$2m$) \ef{QQ1}.
 Quasilinear Schr\"odinger-type models, some of which can be
  written as
    \beq
    \label{Q1}
    \ii u_t + \D u + \b |u|^{p-1}u + \theta (\D |u|^2)u=0 \inA,
 \eeq
    are not a
   novelty in several physical situations such as superfluid theory, dissipative quantum mechanics,
    and    in turbulence
   theory. We refer to some papers from the 1970s and 1980s \cite{Hasse80, Por76, Spat77}, to
   \cite{Guo05} for further references,
    and to \cite{Lange99} for a
   more mathematical knowledge, as
   well as to \cite{Maj97} and to Zakharov's {\em et al} papers \cite{Zak04, Zak01}, as
   a sufficient  source of other reference and deep physical results.
Note that the quasilinear model proposed in \cite{Maj97} in 1997
is more related to a kind of a ``$p$-Laplacian operator" structure
with fractional derivatives, such as
 \beq
 \label{QQ2}
  \tex{
 \ii u_t= \l |D_x|^{\frac \b 4} \big(\big| |D_x|^{\frac \b
 4}u\big|^2 |D_x|^{\frac \b 4} u \big) +|D_x|^\a u,
 }
  \eeq
  with a standard Fourier-definition of operators $|D_x|^\a$,
  having the symbol $|\xi|^\a$, so that, in particular $|D_x|^2 = - D_x^2>0$.
 Then $\l=1$ in \ef{QQ2} corresponds to the original defocusing model. Here, real
parameters $\a$ and $\b$ control dispersion and nonlinearity
respectively. The standard NLSE then occurs for $\a=2$ and $\b =0$
in \ef{QQ2}.

 For $n=0$, the QLSE \ef{QQ1} formally reads as the
linear original one \ef{1}. Indeed, the QSLE \ef{QQ1} still
remains rather rare and seems even an exotic equation in PDE
theory, and therefore we now are going to propose a general
approach to understanding of its new internal
properties.\footnote{While being in Bath in 2008, Peter Markowitz,
answering a question of the first author, fast and witty called
\ef{QQ1}, $m=1$,
 an {\em NLSE of
the porous medium type} bearing in mind the classic parabolic PME
 $$ 
 u_t= \D (|u|^n u) \inA \whereA  n>0.
 $$ 
 }


Overall, bearing in mind some clear discrepancies between the
physical quasilinear models and the proposed one \ef{QQ1}, we just
say that this one was chosen as a typical example only to
demonstrate our branching approach, so that others, more
physically motivated models, would do the same, when the main
mathematical ideas would  have been properly explained and
motivated.

 \subsection{Two ``adjoint" nonlinear eigenvalue problems}

 Thus, as in the linear case for $n=0$, we are going to study
 {\em global asymptotic behaviour} (as $t \to +\iy$)
 and
{\em finite-time blow-up behaviour} (as $t \to T^-<+\iy$)
of solutions of the QLSE \ef{QQ1}.


Overall, we are looking
 for similarity solutions of \ef{QQ1} of  two ``forward" and Sturm's ``backward"  types:

 \ssk

  (i) {\em global
 similarity patterns} for $t \gg 1$, and

 \ssk

 (ii) {\em blow-up
 similarity ones} with the finite-time behaviour as $t \to
 T^-<\iy$.

 \ssk

 Both classes of such particular solutions of the QSLE \ef{QQ1}
 are written in the joint form as follows, by setting $T=0$ in (ii):
  \beq
  \label{upm}
   \tex{
 u_{\pm}(x,t)= (\pm t)^{-\a} f(y), \quad y=x/(\pm t)^\b \whereA \b= \frac {1-\a
 n}{2m},
 \,\,\, \mbox{for} \,\,\,
  \quad \pm t > 0,
  }
  \eeq
  where  similarity profiles $f(y)$ satisfy the following {\em
  nonlinear eigenvalue problems}, respectively,
   \beq
   \label{pm}
    \fbox{$
   \tex{
   {\bf (NEP)_\pm}: \quad
 \BB^\pm_n (\a,f) \equiv
- \ii (-\D)^m (|f|^n f)
 \pm \b y \cdot \n f \pm
 \a f=0 \inB \ren.
 }
  $}
 \eeq
 Here, $\a \in \re$ is  a parameter, which stands in both cases for admitted {\em
 real} (!)
  {\em nonlinear eigenvalues}.
Thus, the sign ``$+$", i.e., $t>0$, corresponds to global
asymptotics as $t \to +\iy$, while ``$-$" ($t<0$) yields blow-up
limits $t \to T=0^-$ describing  a ``micro-scale" structures of
the PDE. In fact, the blow-up patterns are assumed to describe the
structures of ``multiple zeros" of solutions of the QLSE. As we
have mentioned, this idea goes back to Strum's analysis of
solutions of the 1D heat equation performed in 1836 \cite{St}; see
\cite[Ch.~1]{GalGeom} for the whole history and applications of
these fundamental Sturm's ideas and two {\em zero set Theorems}.

 Being equipped with proper ``boundary conditions at infinity",
 namely,
  \beq
  \label{bc1}
  \mbox{for global case}, \quad \BB^+_n(\a,f):
  \quad f(y) \,\,\,\mbox{is ``maximally" oscillatory as $y \to \iy$},
   \quad \mbox{and}
  \eeq
   \beq
  \label{bc2}
  \mbox{for blow-up case}, \quad \BB^-_n(\a,f):
  \quad f(y)\,\,\,\mbox{has a  ``minimal growth" as $y \to \iy$},
  \eeq
equations \ef{pm} produce  {\em true} two {\em nonlinear
eigenvalue problems} to study, which can be considered as a pair
of mutually ``adjoint" ones. Note that \ef{bc2} actually also
means that the admitted nonlinear eigenfunctions {\em are not} of
a type of a maximal oscillatory behaviour at infinity that
connects us with the issue \ef{mm10} (which, however, is not
sufficient, and a growth analysis at infinity should be involved
in parallel, as shown below).

Let us discuss in greater detail the meaning of those above
conditions at infinity. Firstly, \ef{bc1} means that, due to the
type of nonlinearity $|f|^n$, the oscillatory component such as
\ef{mm4} is admissible, with, of course, an extra generated
algebraic factor of the WKBJ-type, which we do not specify hereby.
This can be explained as follows: if $f(y)$ has a standard
WKBJ-type two-scale asymptotics
 \beq
 \label{xx1}
 f(y) \sim |y|^\d {\mathrm e}^{a |y|^\a} \asA y \to \iy,
  \eeq
  then, since $|{\mathrm e}^{a |y|^\a}|=1$ for $a \in \ii \re$,
   substituting into $\ef{pm}_+$ yields  the balance
 \beq
 \label{xx2}
 \d n+(2m-1)(\a-1)=1,
  \eeq
  i.e., different from the purely linear one as in \ef{mm2},
  since, for $n >0$,
  the exponent $\d$ from the slower varying factor is involved (we
  do not calculate it here being a standard asymptotic procedure).

Secondly, \ef{bc2} assumes actually also a ``minimal" growth at
infinity. Namely, quite similar to the linear problem for $n=0$,
the first two terms in $\ef{pm}_-$ generate a fast growing bundle:
as $y \to \iy$ (as usual, we omit slower oscillatory components)
 \beq
 \label{bc21}
- \ii (-\D)^m (|f|^n f)
 - \b y \cdot \n f +...=0  \LongA
 f(y) \sim |y|^{\frac{2m}n}.
  \eeq
  On the other hand, two linear terms in $\ef{pm}_-$ lead to a different slower
  growth as $y \to \iy$:
   \beq
   \label{bc22}
...\, - \b y \cdot \n f -
 \a f=0 \LongA f(y) \sim |y|^{-\frac \a \b} \equiv |y|^{\frac
 {2m|\a|}{1+|\a|n}}
  \eeq
  (recall that $\a(0)=\l <0$). Since
   \beq
   \label{bc23}
    \tex{
\frac
 {2m|\a|}{1+|\a|n}< \frac {2m}n,
  }
  \eeq
 this actually means that \ef{bc2} establishes a kind of a
 ``minimal" growth of admissible nonlinear eigenfunctions at
 infinity corresponding to \ef{bc22}. For $n=0$, this implies a {\em polynomial}
growth, and
 all the admissible (extended) eigenfunctions of ${\mathbb B}^*$ turned out to be
 generalized Hermite polynomials \ef{psi**1}.
Note that,  in self-similar approaches and ODE theory, such
``minimal growth" conditions are known to define similarity
solutions of the {\em second kind}, a term, which was introduced
by Ya.B.~Zel'dovich in 1956
 \cite{Zel56}, and many (but indeed easier) such ODE problems have been rigorously
 solved since that.
  For quasilinear problems such as
 \ef{pm}, the condition \ef{bc2} is incredibly more difficult.
We thus cannot somehow rigorously justify that the problem
$\ef{pm}_-$, \ef{bc2} is well posed and admits a countable family
of solutions and nonlinear eigenvalues $\{\a^-_\g(n)\}$. Actually,
the homotopy deformation as $n \to 0^+$ is the only  our original
intention {\em to avoid} such a difficult ``direct" mathematical
study of this nonlinear blow-up eigenvalue problem.

 All related aspects and notions used
above and remaining unclear will be properly discussed and
specified.

Of course, these conditions \ef{bc1} and \ef{bc2} remind us the
``linear" ones associated with \ef{mm9} for ${\mathbb B}$ and
\ef{mm10} for ${\mathbb B}^*$ respectively, justified earlier for
$n=0$. Indeed, a better understanding of those conditions in the
nonlinear case $n>0$ demands a much more difficult mathematics.
However, one can observe that both \ef{bc1} and \ef{bc2} are just
two {\em asymptotic} (not global ones) problems concerning
admitted
 behaviour of solutions of \ef{pm} as $y \to \iy$, so that, at
 this moment we are in a position to neglect these and to face
 more fundamental issues to be addressed below. Note also that, at
 least, in 1D or for radially symmetric solutions in $\ren$, such asymptotic problems
 for not that hard nonlinear ODEs are {easily} solvable.

Thus, for $n=0$, equations \ef{pm}, equipped with proper weighted
$L^2$ spaces, take the very familiar form: the  corresponding
differential expressions are
 \beq
 \label{QQ6} \BB, \,\,\, \BB^*: \quad
  - \ii (-\D)^m f
 \pm \b y \cdot \n f \pm
 \a f=0 \inB \ren.
  \eeq
Then we observe the obvious relation between $\a$'s and $\l$'s
from the spectrum $\s({\mathbb B})=\s({\mathbb B}^*)$:
 \beq
 \label{QQ7}
  \tex{
 \forA {\mathbb B}: \quad \a= - \l + \frac N{2m} \andA \forA {\mathbb B}^*: \quad
 \a=\l.
 }
 \eeq

\ssk

Thus, our next goal is to show, by using any means, that, at least
for small $n>0$, the nonlinear eigenvalue problems
 \beq
 \label{main1}
 \mbox{
 ${\bf (NEP)_\pm}$ admit  countable sets of
 solutions $\Phi^\pm(n)=\{\a_\g^\pm,f^\pm_\g\}_{|\g| \ge 0}$,}
  \eeq
  where, as usual and as it used to be in the linear case,
   $\g$ is a multiindex in $\ren$ to numerate   the pairs.

The last question to address is whether these sets
 \beq
 \label{main2}
 \mbox{
 $\Phi^\pm(n)$
of nonlinear eigenfunctions are {\em evolutionary complete}},
 \eeq
i.e., describe {\em all} possible asymptotics as $t \to +\iy$ and
$t \to 0^-$ (on the corresponding compact subsets in the variable
$y$ in \ef{upm}) in the CP for the QLSE \ef{QQ1} with bounded
integrable (and possibly compactly supported, -- any assumption is
allowed) initial data. Our main approach is the idea of a
``homotopic deformation" of \ef{QQ1}  as $n \to 0^+$ and reducing
it to our linear equation \ef{1},
 for which both problems \ef{main1} and \ef{main2} are solved
 positively by a non-standard and not self-adjoint spectral
 theory of the linear operator pair $\{{\mathbb B},{\mathbb B}^*\}$.


 \subsection{Example: first explicit nonlinear eigenfunctions}

For $m=1$, the problem $\ef{pm}_+$ has the following first pair
associated with the explicit kernel \ef{Ga1}:
 \beq
 \label{QQ9}
  \tex{
  \a_0^+(n)=\frac N{2+Nn} \andA f_0^+(y)=\big(\frac
  2{2+Nn}\big)^{\frac 1n}\, {\mathrm e}^{\frac{\ii |y|^2}4}.
   }
   \eeq
Not that surprisingly, regardless the degeneracy of the QLSE
\ef{QQ1} at the zero level $\{u=0\}$, \ef{QQ9} shows that the
solution exhibits no {\em finite interfaces}. This is in a
striking difference with, say, as a typical example, the TFE--4
 \beq
 \label{tfe4}
 u_t=-\n \cdot (|f|^n \n \D f) \inA \whereA n>0,
  \eeq
which is known to admit compactly supported solutions in both the
FBP and the Cauchy problem setting; see \cite{GiacOtto08, Grun04}
and, respectively, \cite{EGK2, EGK3}, as a source of main results
and further references.

Concerning the ``blow-up problem",  for any $m \ge 1$, the adjoint
nonlinear eigenvalue one $\ef{pm}_-$ has the obvious first pair
 \beq
 \label{QQ10}
  \tex{
  \a_0^-(n)=0 \andA f_0^-(y) \equiv 1,
  }
  \eeq
  i.e., the same as for $n=0$, where $\psi_0^*(y) \equiv 1$, as
  the first Hermite polynomial; see \ef{psi**1}.

Anyway, in our further analysis, we cannot rely on any explicit
representation of any nonlinear eigenfunctions. So, we now very
briefly begin to explain our approach.

\subsection{Branching for the ``forward" problem for $t \gg 1$}

We perform our bifurcation-branching analysis following the lines
of classic theory \cite{Berger, Deim, KrasZ, VainbergTr}, etc., in
the case of {\em finite} regularity. However, we must admit that
this classic one does not cover rigorously  the type of
$n$-branching as $n \to 0$, which requires extra difficult study,
which we cannot address here.

Thus, we are looking for a countable set of nonlinear pairs
$\{a_\g^+(n),f^+(y;n)\}_{|\g| \ge 0}$, which are assumed to
describe {\em all} (hopefully, for $n>0$ small) asymptotic
patterns for the QLSE \ef{QQ1} as $t \to +\iy$, {\em up to
possible centre manifold patterns}; see Section \ref{S8.2}.
However, as we have seen, such special patterns do occur, if there
is a certain {\em interaction} (a {\em transitional behaviour})
between the linear and the nonlinear terms in equations such as
\ef{eqPer}. This looks rather unreal for the current problem under
consideration. However, the {\em evolutionary completeness} of the
nonlinear patterns $\Phi^+(n)$ as $t \to +\iy$ remains an
extremely difficult open problem, which probably will  be
extremely hard to solve completely rigorously.

Thus, assuming that $n>0$ is small, we perform asymptotic
expansions in the operators and coefficients in $\ef{pm}_+$.
Evidently, the crucial one is in the nonlinearity, which requires:
 \beq
 \label{u1}
 |f|^n f \equiv f{\mathrm e}^{n \ln |f|}= f(1+n \ln |f| + o(n))
 \asA n \to 0^+.
  \eeq
It is clear that the neighbourhoods of the nodal set of $f(y)$ are
key for \ef{u1} to be valid in any weak sense (precisely this is
needed for the equivalent analysis of the inverse {\em integral}
compact operators involved, where a proper justification must take
place). If $f(y)$ has a {\em nice nodal set consisting of a.a.
isolated and ``transversal" a.e. zero surfaces $($or just some
points only$)$, with no ``thin" concentration subsets}, the
expansion \ef{u1} can be valid even in the standard pointwise
sense, or at least in the weak sense. However, we do not know and
currently cannot prove such deep properties of the nonlinear
eigenfunctions involved. Note that, for $n=0$, the generating
formula \ef{eigen}, with a proper knowledge of such nice zero set
properties of the rescaled fundamental kernel $F(y)$ (this is
doable), guarantees
 such necessary properties of eigenfunctions $\psi_\g(y)$.

 The above discussion establishes the main
 hypothesis to make our branching analysis to be (almost)
 rigorous, being applied, of course, to the equivalent integral
equation, where establishing  some further compact and other
necessary properties of the nonlinear integral operators would
take some time, indeed. This can be also done, since the problems
$\ef{pm}_\pm$ can be  reduced to a semilinear form.

The rest of the expansions in $\ef{pm}_+$ are straightforward
(here, we already fix by $\g$ an $n$-branch we are going to trace
out):
 \beq
 \label{u2}
  \begin{matrix}
  \a_\g(n)=\a_\g(0)+ \hat \a_\g n +... \equiv -\l_\g+ \frac N{2m}+
  \hat \a_\g n+...\, , \ssk\ssk \\
   \frac {1-\a_\g(n)n}{2m}= \frac 1{2m} - a_\g n+...\, , \,\,\,
   \mbox{where}
   \,\,\,a_\g= \frac {-\l_\g+ \frac N{2m}}{2m},
    \end{matrix}
 \eeq
 where we have already omitted all $o(n)$-terms
 (this again assumes extra regularity hypothesis already discussed above).
In the first line in \ef{u2}, the parameter $\hat \a_\g$ is an
extra unknown.

Substituting all the expansions into $\ef{pm}_+$ yields the
following perturbed problem:
 \beq
 \label{u3}
  \tex{
  ({\mathbb B}- \l_\g I)f+ n \, h+...=0, \,\,\, \mbox{where}
  \,\,\, h=  \big[\ii (-1)^{m+1} \D^m(f \ln|f|)-a_\g y
  \cdot \n f + \hat \a_\g f\big].
   }
    \eeq
Thus, as $n \to 0^+$, we must look for a solution close to the
eigenspace
 \beq
 \label{u4}
  \tex{
 {\rm ker} \,({\mathbb B}-\l_\g I)={\rm Span}\,\{
 \psi_\g,\,\,|\g|=l\}.
 }
 \eeq
 Therefore, under prescribed hypothesis, solutions take the form
 \beq
 \label{u41}
 \tex{
 f= \sum_{|\s|=l} c_\s \psi_\s + n \phi_\g+...\,,
 }
 \eeq
where the expansions coefficients $\{c_\s\}_{|\s|=l}$ and the
orthogonal part $ \phi_\g$ are unknowns.

Finally, substituting \ef{u41} into \ef{u3} yields the
$O(n)$-problem
 \beq
 \label{u42}
 ({\mathbb B}-\l_\g I) \phi_\g + h=0.
  \eeq
  Thus,
 the necessary (and, in properly regular cases, the sufficient)  orthogonality condition
 of the solvability of \ef{u42},
 \beq
 \label{u5}
  h \,\,\bot^* \,\,{\rm ker} \,({\mathbb B}-\l_\g I),
   \eeq
yields the following Lyapunov--Schmidt scalar branching equation:
for all $|\d|=l$,
 \beq
 \label{u6}
  \tex{
   \langle \ii (-1)^{m+1} \D^m(\sum c_\s \psi_\s \ln|\sum c_\s \psi_\s|)-a_\g y
  \cdot \n \sum c_\s \psi_\s + \hat \a_\g \sum c_\s \psi_\s, \,
  \psi_\d^* \rangle_* =0,
    }
     \eeq
     where an integration by parts in the second term can be
     performed to simplify the expressions; see the next Appendix
     D for further details concerning the meaning of such extended
     generalized linear functionals. We must admit again that a proper
     well posed way for calculating values of such extended functionals
     is not available, so we present such a branching analysis  just as
     an example reminding analogies with the classic approaches.

Overall, \ef{u6} is the required {\em algebraic system} for the
unknowns $\hat \a_\g$ and $\{c_\s\}$ (a convenient normalization
condition on the latter expansion coefficients may be added). As
usual, once the branching equation  \ef{u6} has been properly
solved, this allows one  to get the corresponding unique solution
of the differential equation \ef{u42}, etc.

Indeed, \ef{u6} is a very difficult algebraic system, which is not
of any variational form, so one cannot use powerful category-genus
theory \cite{Berger, KrasZ} to predict a number of solutions,
i.e., a number of such $n$-branches originated from the given
eigenspace in \ef{u4}. Note that, for $m=1$, when there exists
some ``symmetry" of differential forms ${\mathbf B}$, ${\mathbf
B}^*$ reflected in \ef{ss1} and \ef{kam1}, the system \ef{u6}
reveals some ``variational-like" features, since the
eigenfunctions $\psi_\g$ and polynomials $\psi_\g^*$ can be
identified in a weighted space. However, this is supposed to
happen in a space with indefinite metric, so we do not check how
this can be helpful.

It is worth mentioning that the problem on a sharp estimate of a
number of  solution branches emanating from  an eigenspace,
remains essentially open even for classes of well-understood
variational operators. On one hand, the first conclusion is
classic: the number of branches is not less than the dimension of
the eigenspace: indeed, since the corresponding algebraic system
(like \ef{u6}) remains also variational, the category of the
functional set is not less than the linear eigenspace dimension,
whence the result. But, obviously, a sharper estimate of the
solutions number becomes essentially nonlinearity-dependent, so
that
 this is not (and, possibly,  cannot be  in the maximal generality)
completely understood. For  a number of branches that can emanate
from the trivial solution,   there have been  obtained  some
specific examples only. There exist some results for potential
operators (see \cite{DGM} and \cite{R} as a guide), and a very few
for non-gradient and non-self-adjoint operators \cite{KHK}. For
the TFE--4 \ef{tfe4}, such a branching analysis as $n \to 0^+$
reveals a lot of technical difficulties, though some problems for
simple and semisimple eigenvalues are shown to admit a rather
definite path towards, \cite{Pablo2}.

\subsection{Blow-up scaling problem: $t \to 0^-$}

This is quite similar. Since the nonlinear eigenvalue equations
$\ef{pm}_\pm$ differ by the sign in the linear terms only (but the
boundary-radiation conditions at infinity are entirely different),
instead of \ef{u3}, we arrive at
 \beq
 \label{u3-}
  \tex{
  (\BB^*- \l_\g I)f+ n \, h+...=0, \,\,\, \mbox{where}
  \,\,\, h=  \big[\ii (-1)^{m+1} \D^m(f \ln|f|)+a_\g y
  \cdot \n f - \hat \a_\g f\big].
   }
    \eeq
Therefore, looking for solutions \ef{u41}, now over the kernel of
the adjoint operator $({\mathbb B}^*-\l_\g I)$,
 \beq
 \label{u41-}
 \tex{
 f= \sum_{|\s|=l} c_\s \psi_\s^* + n \phi_\g^*+...\,,
 }
 \eeq
 we arrive at the ``adjoint" algebraic system:
for all $|\d|=l$,
 \beq
 \label{u6-}
  \tex{
   \langle \ii (-1)^{m+1} \D^m(\sum c_\s \psi_\s^* \ln|\sum c_\s \psi_\s^*|)+a_\g y
  \cdot \n \sum c_\s \psi_\s^* - \hat \a_\g \sum c_\s \psi_\s^*, \,
  \psi_\d \rangle_* =0,
    }
     \eeq
which are not easier that the first one \ef{u6}, and does not have
any variational structure, so the same principal difficulties on
the solvability (this is easier to do in the 1D or the radial
case) occurs, and especially on the number of solutions.


Recall that the patterns
 \beq
 \label{u7-}
  \tex{
u(x,t) \sim {\mathrm e}^{\a_\g(n) \t} \sum_{|\s|=l} c_\s
\psi_\s^*(y)+...\, \whereA y= \frac x{(-t)^{\b_\g(n)}}, \quad \t =
- \ln(-t),
 }
 \eeq
 are assumed to describe, for all finite multiindices $\g$ and all
 admitted solutions of the branching equation \ef{u6-}, the whole
 variety of ``micro-scale patterns", which are available for the
 QLSE \ef{QQ1} at a given arbitrary point $(0,0^-)$, unless  some
 centre-subspace-type patterns might appear due to nonlinearities
 involved. In any case, we believe that \ef{u7-} describe (at
 least, a.a.)
 of generic formations of ``multiple zeros" of solutions. In other
 words, for small $n>0$, multiple zeros of the ${\rm Re}\, u(x,t)$ and/or ${\rm
 Im}\, u(x,t)$ as $(x,t) \to (0,0^-)$ are created by a
 self-focusing of zero surfaces of the corresponding profiles
 $f(y)$ given in \ef{u41-}, which comprises a proper linear
 combination of the generalized Hermite polynomials \ef{psi**1}.
 Indeed, this study requires further extensions, however, a
 complete and fully rigorous answer seems cannot  be achieved.

 \ssk

With such an ``optimistic point", we end up the
 list of possible
applications of our refined scattering linear spectral $\{{\mathbb
B},{\mathbb B}^*\}$-theory for $2m$th-order rescaled
Schr\"odinger operators.

\end{small}
\end{appendix}

\begin{appendix}
\begin{small}

\section*{Appendix D. Eigenfunction expansions and little Hilbert spaces}
 \label{SHH}
 \setcounter{section}{4}
\setcounter{equation}{0}


 Given below further developing of spectral theory  is not
 necessary for our main applications concerning classification of all the global and blow-up
asymptotics, which follow from the corresponding spectral
decompositions of semigroup representation of solutions. However,
we think that these accompanying results are interesting and
actually allow to extend the results to wider classes of
solutions.

Recall that the main difficulty with a proper definition of the
operator ${\mathbb B}$ \ef{BB11} by its spectral decomposition
\ef{de1} was associated with the fact that the necessary for us
eigenfunctions $\{\psi_\b\}$ were {\em extended}, i.e., did not
belong to the present domain. In other words, these
eigenfunctions, which inevitably appeared in eigenfunctions
expansions \ef{w1New} of rather ``good" solutions of Schr\"odinger
equations, were much ``worse" that the solutions themselves. To
get rid of such a controversy and to restore the true meaning of
the operator pair $\{\BB,\, \BB^*\}$ (instead of its restriction
$\{{\mathbb B}, {\mathbb B}^*\}$), we first introduce new {\em
extended spaces of closures}.

\subsection{Subspace where  $\Phi $ is closed}



Given the complete subset $\Phi =\{\psi_\b\}$ for the non
self-adjoint operator ${\bf B}$, we define the  linear subspace
$\tilde L^2_\rho$ of eigenfunction expansions,
 \begin{equation}
\label{UU1}
 \tex{
 v \in \tilde L^2_\rho \quad \mbox{iff} \quad v = \sum
c_\b \psi_\b \,\,\,  \mbox{with convergence in} \,\,\,
L^2_\rho(\ren),
 }
 \end{equation}
as the closure of the subset of finite sums
 \beq
 \label{fs1}
 \tex{
\big\{\sum_{|\b| \le K} c_\b \psi_\b, \, K \ge 0\big\}
 }
 \eeq
 in the
$L^2_\rho$-norm.






For clarifying such a space, we now derive  better estimates to
see which $\{c_\b\}$ satisfy \ef{UU1}. Namely, we will use the
equality \ef{mm6} for the key exponent in the exponential
representation \ef{mm4} of the asymptotics (which is responsible
for the sharp estimate \ef{ff1} of eigenfunctions $\psi_\b$). Then
similar to \ef{jj3}, but sharper,  we then obtain for $l=|\b| \gg
1$,
 \beq
 \label{es60}
  \tex{
  \int \rho |\psi_\b|^2 \sim \big( \frac l {\mathrm e}\big)^{-l} l^{\frac
  {2l(\a-1)}\a} \big[ \frac{2(\a-1)}{\a l} \big]^{\frac
  {2l(\a-1)}\a} (2m)^{-\frac {2l}{2m-1}}
  = l^{- \frac {l(2-\a)}\a} \big[ \frac{ {\mathrm e}^{2m-1}}{2m(m
   {\mathrm e})^{(2m-1)/m}}\big]^{\frac l{2m-1}}.
   }
   \eeq
 Next, since
  \beq
  \label{es61}
   \tex{
    \|v\|^2_{L^2_\rho} = \sum\limits_{(\b,\g)} c_\b c_\g \int\limits_{\ren} \rho
    \psi_\b\psi_\g,
     }
     \eeq
  \ef{es60} implies that:
   \beq
   \label{es62}
    \tex{
    v \in \tilde L^2_\rho \,\,\,\mbox{if\,\, $c_\b$ does not glow faster
than\,\, $l^{2(\frac{2-\a}\a-\e)}$ for $l=|\b| \gg 1$},
 }
 \eeq
 where $\e>0$ is any arbitrarily small constant. Of course,
 \ef{es60} defines more optimal and weaker inclusion conditions,
 but \ef{es62} clearly explains how this works.


\subsection{Bi-orthonormality of the bases}
 Obviously, this is a principal issue for all the applications, where eigenfunction expansion
  techniques
 take part.
 As we have seen, the eigenfunctions expansions such as
 \ef{BB12} introduce standard linear functionals $\langle w_0,
 \psi_\b^* \rangle$,  which are well defined for all functions $w_0
 \in \LLL$, so that, as usual,  $\psi_\b^*$ is an element of the adjoint space
 $\LL$, with $\rho= \frac 1 {\rho^*}$, as customary.


\ssk

As the next step, according to our construction above, we have to
define some generalized {\em extended linear functionals}
 from the adjoint
space $\tilde L^{2*}_\rho$. On one hand,
 this would correspond
 to a standard
procedure of extension of such continuous uniformly convex
functionals
 by the
Hahn--Banach classic theorem in linear normed
spaces,
 \cite{KolF}. As we have seen earlier,
those linear functionals are well defined according to \ef{MBB1}
in $L^2_{\rho^*}$ with the standard (not in any {\em v.p.} or a
canonical regularized, etc., sense) definition of the integrals.

On the other hand, such extended linear functionals cannot be
understood in a standard sense, so that we refer to them as to
{\em generalized} ones. It seems, a full proper definition of such
extended linear functionals in a usual functional framework will
require a  deeper analysis of the actual functional spaces and
metric/topologies involved, which will essentially decline us from
main PDE applications [Especially, since in some of the
applications, we do not and even cannot pretend to be
mathematically rigorous.]




 Therefore, in other words, for any $v \in \tilde L_\rho^2$,
we define  {\em extended linear  functionals} for any $\b$ as:
 \beq
 \label{MEAN1}
 \mbox{
 $\langle v, \psi_\b^* \rangle_* \equiv \langle
v, \bar\psi_\b^* \rangle$\, denotes  the expansion coefficient
$c_\b$ of $v$ in (\ref{UU1}).
 }
 \eeq
 In view of the performed construction of the space
$\tilde L_\rho^2$ via closure of finite sums (\ref{fs1}), it is
not difficult to see that such generalized continuous linear
functionals are defined uniquely (in view of the density of finite
sums (\ref{fs1})).

Overall, in the sense of \ef{MEAN1},
 the standard bi-orthonormality of the bases $\{\psi_\b\}$ and $\{\psi^*_\b\}$ becomes trivial:
 \beq
 \label{bi1}
 \langle \psi_\b,  \psi_\g^* \rangle_* \equiv
\langle \psi_\b, \bar \psi_\g^* \rangle = \d_{\b\g} \quad
\mbox{for any} \,\,\, \b \,\, \mbox{and} \,\, \g,
 \eeq
 where $\langle \cdot,\cdot \rangle$  is the usual duality
product in $L^2(\ren)$ and $\d_{\b\g}$ is the Kronecker delta.


 Similarly, using the subset $\Phi ^*=\{\psi_\b^*\}$
of the generalized Hermite polynomials \ef{psi**1},
  we are obliged to define
 the corresponding subspace $\tilde L^{2}_{\rho,*}$ of eigenfunction expansions, and eventually treat
 similarly
   the adjoint extended linear functionals
$\langle w, \psi_\b \rangle_*$ for any $w \in \tilde
L_{\rho,*}^{2}$.



In writing \ef{bi1}, we use the standard $L^2$-metric, which is
convenient to see the normalization since $\psi_\b$ is essentially
the $D^\b$ derivative and $\psi_\b^*$ is a polynomial, so that,
including normalization factors yields
 \beq
 \label{bi2}
  \tex{
\langle \psi_\b, \bar \psi_\b^* \rangle = \frac
{(-1)^{|\b|}}{\b!}\,
 \int\limits_{\ren} D^\b F(y) \,\, (y^\b+...)\, {\mathrm d}y= 1
 \quad \mbox{for any multiindex} \,\,\, \b
 }
 \eeq
 in the sense of formal integration by parts. Therefore, the integral
 itself is not of a standard meaning, but can be  treated in an involved ``{\em v.p.}-like sense", which
 is difficult to clarify, and we do not feel any actual necessity to do this.
 A similar formalism exists for the whole  eigenfunction set occurred in  \ef{bi1}.
 Thus,
 we will use \ef{bi1} in the eigenfunction expansions to follow,
  bearing in mind its  actual meaning specified above in (\ref{MEAN1}).


\subsection{Little Hilbert and Sobolev spaces}

It is convenient to introduce a {\em little}  Hilbert space
$l^2_\rho$ of functions $v = \sum a_\b \psi_\b \in \tilde
L^2_\rho$ with coefficients satisfying
 \begin{equation}
 \label{abeta1}
  \tex{
 \sum |a_\b|^2 < \infty,
  }
 \end{equation}
where the scalar product  and the induced norm are given by
 \begin{equation}
 \label{vwin}
  \tex{
  (v,w)_{0} = \sum a_\b \bar c_\b, \,\, w = \sum c_\b \psi_\b \in l^2_\rho
  \quad \mbox{and} \quad
\|v\|^2_0 = (v,v)_0.
 }
 \end{equation}
 Therefore, $\Phi $ is now treated as
 a Riesz basis in $\tilde L^2_\rho(\ren)$, \cite{BS, GGK}.
We next define a little Sobolev space $h^{2m}_\rho$ of functions
$v \in l^2_\rho$ such that
 $$
  \tex{
  \sum |\l_\b c_\b|^2 < \infty.
  }
  $$

  The scalar product and the induced norm in $h^{2m}_\rho$ are
   \begin{equation}
   \label{5656}
    \tex{
  (v,w)_{1} = (v,w)_0 + ({\bf B} v,{\bf B} w)_0, \quad \|v\|_{1}^2 =
  (v,v)_{1} \equiv \sum \big(1+|\l_\b|^{2}\big)|c_\b|^2,
  }
  \end{equation}
  where our bounded operator $\BB: h^{2m}_\rho \to l^2_\rho$ has the
meaning $\BB: \{c_\b\} \to \{\l_\b c_\b\}$. This norm is
equivalent to the graph norm induced by the positive operator
$(-{\bf B} + a I)$ with  $a>0$. Then $h^{2m}_\rho$ is the domain
of ${\bf B}$ in $l^2_\rho$, and, by Sobolev's embedding theorem,
 \begin{equation}
  \label{SobE}
   \tex{
 h^{2m}_\rho \subset l^2_\rho \quad \mbox{compactly},
 }
 \end{equation}
 which follows from the criterion of compactness in $l^p$,
 \cite{KolF}.



\subsection{Basic properties in $l^2_\rho$}

Firstly, it follows that  ${\bf B}$  {\em is self-adjoint}
(symmetric) in $l^2_\rho$,
 \begin{equation}
 \label{Bsad}
 ({\bf B} v, w )_0 =  (  v, {\bf B} w )_0 \quad
 \mbox{for all} \,\,\, v,w \in h^{2m}_\rho.
 \end{equation}
 Secondly,
 we state some other
 straightforward
 consequences.

\begin{proposition}
\label{PrDens} {\rm (i)} The Hilbert space  $l^2_\rho$ is a dense
subspace of $\tilde L^2_\rho$ in $L^2_\rho(\ren)$;

{\rm (ii)} $\Phi =\{\psi_\b\}$ is complete and closed in
$l^2_\rho$ in the topology of $L^2_\rho(\ren)$;

{\rm (iii)}  the resolvent  $({\bf B} - \l I)^{-1}$ for $\l \not
\in \s({\bf B})$ is compact in $l^2_\rho$; and

{\rm (iv)} ${\bf B}$ is sectorial in $l^2_\rho$.
\end{proposition}

\noi {\em Proof.} (i) $l^2_\rho$ is separable and complete since
the same is true for the isomorphic Hilbert space $l^2$ of
sequences. Let us show that $l^2_\rho\subseteq \tilde L^2_\rho$.
For any $v
 \in l^2_\rho$,
  $$
   \tex{
 \int \rho |v|^2 \, {\mathrm d}y  = \int |\sum a_\b \sqrt \rho \psi_\b|^2 \, {\mathrm d}y
 \le \sum |a_\b|^2 \sum \int {\mathrm e}^{-|y|^\a} \big| \frac 1{\sqrt{\b !}} \, D^\b
 F(y)\big|^2 \, {\mathrm d}y ,
 }
 $$
 and by the same  estimates
as in (\ref{jj1}), \ef{jj2}, we conclude that, for  $l = |\b| \gg
1$,
 \begin{equation}
 \label{est22}
  \tex{
 \int {\mathrm e}^{-|y|^\a}  \big| \frac 1 {\sqrt{\b !}}\, D^\b
 F(y) \big|^2 \, {\mathrm d}y  \le l^{l(-\nu +\e)},
  }
 \end{equation}


 \noi where $\nu = \frac {2-\a}\a  >0$ and $\e > 0$ can be  an arbitrarily small constant.
 Therefore,  $l^2_\rho \subset \tilde L^2_\rho$. Concerning the density of
 $l^2_\rho$,
  we note that
  given a $v= \sum a_\b \psi_\b \in
 \tilde L^2_\rho$, the sequence of truncations $ \{ \sum_{|\b| \le
 K} a_\b \psi_\b, K = 1,2,...\} \subset l^2_\rho$ converges to $v$ in the topology of $L^2_\rho(\ren)$
 as $K \to \infty$
by completeness and closure  of
 $\{\psi_\b\}$.

 (ii) Since $\Phi $ is orthonormal in $l^2_\rho$, it follows
 that the only element orthogonal to $\{\psi_\b\}$ is $0$, and
 hence completeness of $\{\psi_\b\}$ in $l^2_\rho$  follows from the Riesz--Fischer theorem.
  It is closed as an orthonormal subset in a separable
 Hilbert space \cite{KolF}.

(iii) For any $v= \sum a_\b \psi_\b \in l^2_\rho$ from the unit
ball $T_1$ in $l^2_\rho$ with  $\sum |a_\b|^2 \le 1$,
 \begin{align}
 &
  \tex{
 ({\bf B}-\l I)^{-1}v = \sum b_\b \psi_\b,
 } \quad \mbox{where}
 \notag \\
  \tex{
b_\b = \frac {a_\b}{\l_\b-\l}=} - & \, \tex{\frac {a_\b}{|\b|/2m +
\l} = - \frac{2m a_\b}{|\b|}\big[1+
  O(\frac 1{|\b|})\big] \quad \mbox{for} \,\,\, |\b| \gg 1. }  \label{BBB11}
  \end{align}
  Therefore, for any $\e > 0$, there exists $K=K(\e)>0$ such that
  for any $v \in T_1$,
 $$
  \tex{
 \sum\limits _{|\b| \ge K}|b_\b|^2 \le  \frac {4m^2}{K^2} \sum |a_\b|^2
 \le  \frac {4m^2}{K^2} < \e.
  }
 $$
 By the  compactness criterion  in
$l^2$ \cite{KolF},
 $({\bf B} - \l I)^{-1}$ maps $T_1$ onto a compact subset in $l^2_\rho$.

(iv)
 Recall that
 $({\bf B}-\l I)^{-1}$ is a meromorphic function having a
pole $\sim \frac 1\l$ as $\l \to 0$ since $\l_0=0$ has
multiplicity one \cite{GGK}. We then  need an extra estimate on
the resolvent, which is easy to get in $l^2_\rho$ (one can check
that it is not that easy in the big space $L^2_\rho$). In the
sector $\O_\theta= \{\l \in {\mathbb C}:\,\, \l \not =0,
 \,\,|{\rm arg}\,\l| < \frac \pi 2+\theta\}$ with a $ \theta \in
 (0,\frac \pi 2)$, for any
 $v= \sum a_\b \psi_\b \in
l^2_\rho$, we apply (\ref{BBB11}) by using the fact that
 $
  \frac 1{|\l_\b-\l|} \le   \frac 1{|\l|}\,{\sin\theta}$
 in $\O_\theta$ to get
 $$
  \tex{
  \|({\bf B}-\l I)^{-1}v\|_0=  \big(\sum |a_\b|^2 \frac
1{|\l_\b-\l|^2}\big)^{\frac 12} \le \frac 1{\sin\theta} \frac
1{|\l|} \|v\|_0.
 }
 $$
  Since ${\bf B}$ is closed and densely defined, it
 is a sectorial operator in $l^2_{\rho}$, see \cite{Fr}. $\qed$

\ssk

Similarly to (\ref{UU1}), for the adjoint operator ${\bf B}^*$ we
define the subspace $\tilde L^2_{\rho,*} \subset
L^2_{\rho}(\ren)$, where the eigenfunction subset  $\Phi^* $ is
closed, $l^2_{\rho,*}$, $h^{2m}_{\rho,*}$, etc., and next continue
develop similar theory of self-adjointness and other properties



\end{small}
\end{appendix}

\end{document}